\documentclass{gtmon_a}
\pdfoutput=1
\usepackage{pinlabel}

\proceedingstitle{Groups, homotopy and configuration spaces (Tokyo
  2005)}
\conferencestart{5 July 2005}
\conferenceend{11 July 2005}
\conferencename{Groups, homotopy and configuration spaces,
                in honour of Fred Cohen's 60th birthday}
\conferencelocation{University of Tokyo, Japan}

\editor{Norio Iwase}
\givenname{Norio}
\surname{Iwase}

\editor{Toshitake Kohno}
\givenname{Toshitake}
\surname{Kohno}

\editor{Ran Levi}
\givenname{Ran}
\surname{Levi}

\editor{Dai Tamaki}
\givenname{Dai}
\surname{Tamaki}

\editor{Jie Wu}
\givenname{Jie}
\surname{Wu}

\title{Energy of knots and the infinitesimal cross ratio}                    
\author{Jun O'Hara}                  
\givenname{Jun}
\surname{O'Hara}
\address{Department of Mathematics\\Tokyo Metropolitan University\\\newline
Tokyo\\Japan} 
\email{ohara@tmu.ac.jp}                     
\urladdr{  }                       

\volumenumber{13}
\issuenumber{}
\publicationyear{2008}
\papernumber{19}
\startpage{421}
\endpage{445}

\doi{}
\MR{}
\Zbl{}

\arxivreference{0708.2948}


\received{30 May 2006}
\revised{9 May 2007}
\accepted{16 May 2007}
\published{19 March 2008}
\publishedonline{19 March 2008}
\proposed{}
\seconded{}
\corresponding{}
\version{}

\makeatletter
\def\cnewtheorem#1[#2]#3{\newtheorem{#1}{#3}[section]
\expandafter\let\csname c@#1\endcsname\c@thm}

\let\xysavmatrix\xymatrix
\def\xymatrix{\disablesubscriptcorrection\xysavmatrix}
\AtBeginDocument{\let\bar\wbar\let\tilde\wtilde\let\hat\what}

\newtheorem{thm}{Theorem}[section]  
\cnewtheorem{lem}[thm]{Lemma}        
\cnewtheorem{proposition}[thm]{Proposition}
\theoremstyle{definition}
\cnewtheorem{conjecture}[thm]{Conjecture}
\cnewtheorem{wconjecture}[thm]{Willmore Conjecture}
\cnewtheorem{defn}[thm]{Definition}  
\newtheorem*{remark}{Remark}           
\cnewtheorem{problem}[thm]{Problem}

\makeatother  
\makeautorefname{subsection}{Subsection}

\def\a{\alpha}
\def\e{\varepsilon}
\def\w{\wedge}
\def\p{\prime}
\AtBeginDocument{\def\R{\Re\mathfrak{e} \,}}
\def\I{\Im\mathfrak{m} \,}
\def\omegasub#1{\mbox{\large $\omega$}\!\!\mbox{\small${\mbox{\large ${}$}}_{#1}$}}
\def\vect#1{\mbox{\boldmath $#1$}} 
\def\spbmapright#1#2{\smash{%
 \mathop{\hbox to 0.8cm{\rightarrowfill}}
  \limits^{#1}_{#2}}}

\begin{document}

\begin{htmlabstract}
This is a survey article on two topics.  The Energy E of knots can
be obtained by generalizing an electrostatic energy of charged knots
in order to produce optimal knots.  It turns out to be invariant under
M&ouml;bius transformations.  We show that it can be expressed in terms
of the infinitesimal cross ratio, which is a conformal invariant of a
pair of 1&ndash;jets, and give two kinds of interpretations of the real
part of the infinitesimal cross ratio.
\end{htmlabstract}

\begin{abstract}  
This is a survey article on two topics.  The Energy $E$ of knots can
be obtained by generalizing an electrostatic energy of charged knots
in order to produce optimal knots.  It turns out to be invariant under
M\"obius transformations.  We show that it can be expressed in terms
of the infinitesimal cross ratio, which is a conformal invariant of a
pair of $1$--jets, and give two kinds of interpretations of the real
part of the infinitesimal cross ratio.
\end{abstract}

\begin{asciiabstract}  
This is a survey article on two topics.  The Energy E of knots can be
obtained by generalizing an electrostatic energy of charged knots in
order to produce optimal knots.  It turns out to be invariant under
Moebius transformations.  We show that it can be expressed in terms of
the infinitesimal cross ratio, which is a conformal invariant of a
pair of 1-jets, and give two kinds of interpretations of the real part
of the infinitesimal cross ratio.
\end{asciiabstract}

\maketitle

\section{Introduction} 
This is a survey article on two topics, the energy $E$ of knots and
the infinitesimal cross ratio which can give a conformal geometric
interpretation of the energy.

In the first part of this paper we give an introduction to the theory
of energy of knots.  Energy of knots is a functional on the space of
knots which blows up as a knot degenerates to a singular knot with
double points.  It was introduced to produce optimal knots.  The first
example, the energy $E$, was obtained by the author by generalizing an
electrostatic energy of charged knots \cite{OH1}.  Later on,
it was proved to be invariant under M\"obius transformations
(Freedman, He and Wang \cite{Fr-He-Wa}).

The second part of this paper is a survey and an announcement of a
part of the joint work with R\'emi Langevin \cite{LO1,LO2}.  We give a
new interpretation from a viewpoint of conformal geometry.  The
infinitesimal cross ratio is the cross ratio of $x, x+dx, y$, and
$y+dy$, where these four points are considered complex numbers by
identifying a sphere through them with the Riemann sphere $\mathbb
C\cup\{\infty\}$.  It can be considered a complex valued $2$--form on
$K\times K\setminus\Delta$.  It is the unique conformal invariant of a
pair of $1$--jets of a given curve up to multiplication by a constant.
We show that the energy $E$ can be expressed as the integration of the
difference of the absolute value and the real part of the
infinitesimal cross ratio.  We then show that the real part of the
infinitesimal cross ration can be interpreted in two ways: as the
canonical symplectic form of the cotangent bundle of $S^3$ and as a
signed area element with respect to the pseudo-Riemannian structure of
the set of oriented $0$--spheres in $S^3$.

\section{Energy of knots}

\subsection{Motivation}
Just like a minimal surface is modeled on the ``optimal surface'' of a soap film with a given boundary curve, one can ask whether we can define an ``{\em optimal knot}'', a beautiful knot which represents its knot type. 
The notion of energy of knots was introduced for this purpose. 
The basic philosophy is as follows. 

Suppose there is a non-conductive knotted string which is charged uniformly in a non-conductive viscous fluid. 
Then it might evolve itself to decrease its electrostatic energy  without intersecting itself because of Coulomb's repulsive force  until it comes to a critical point of the energy. 
Then we might be able to define an ``optimal knot'' by an embedding that attains the minimum energy within its isotopy class. 
Thus our motivational problem, which was proposed by Fukuhara and Sakuma independently, can be stated as: 

\begin{problem}\label{motivation_probl}\rm (Fukuhara \cite{Fu} and Sakuma
  \cite{Sak})\qua Give a functional $e$ (which we will call an energy) on
  the space of knots $\mathcal{K}$ which satisfies the following
  conditions:
\begin{enumerate}
\item Let $[K]$ denote an isotopy class which contains a knot $K$. 
Define the energy of an isotopy class by $\displaystyle e([K])=\inf_{K^{\p}\in[K]}e(K^{\p})$. 
\item If a knot $K_0$ attains the minimum value of the functional $e$ within its isotopy class, ie, if $e(K_0)=e([K_0]),$ we call $K_0$ an {\em $e$--minimizer} of the isotopy class $[K_0]$. 
\item There is an $e$--minimizer in each isotopy class. 
\end{enumerate}
\end{problem}

\begin{figure}[ht!]
\begin{center}
\labellist\small
\pinlabel $e$ [b] at 90 660
\pinlabel {isotopy class $[K]$} [t] at 175 50
\pinlabel $\{\textrm{immersions}\}$ <20pt, 0pt> [t] at 467 42
\pinlabel* $C^2$--topology [b] at 596 86
\pinlabel* {$\displaystyle \left\{\!\!\begin{array}{l}
            \textrm{singular knots}\\
            \textrm{with double points}
             \end{array}\!\!\right\}$} <10pt, 0pt> [b] at 412 251
\pinlabel (a) [t] <0pt, -15pt> at 311 18
\pinlabel (b) [t] <180pt, -15pt> at 311 18
\endlabellist
\includegraphics[width=0.48\linewidth]{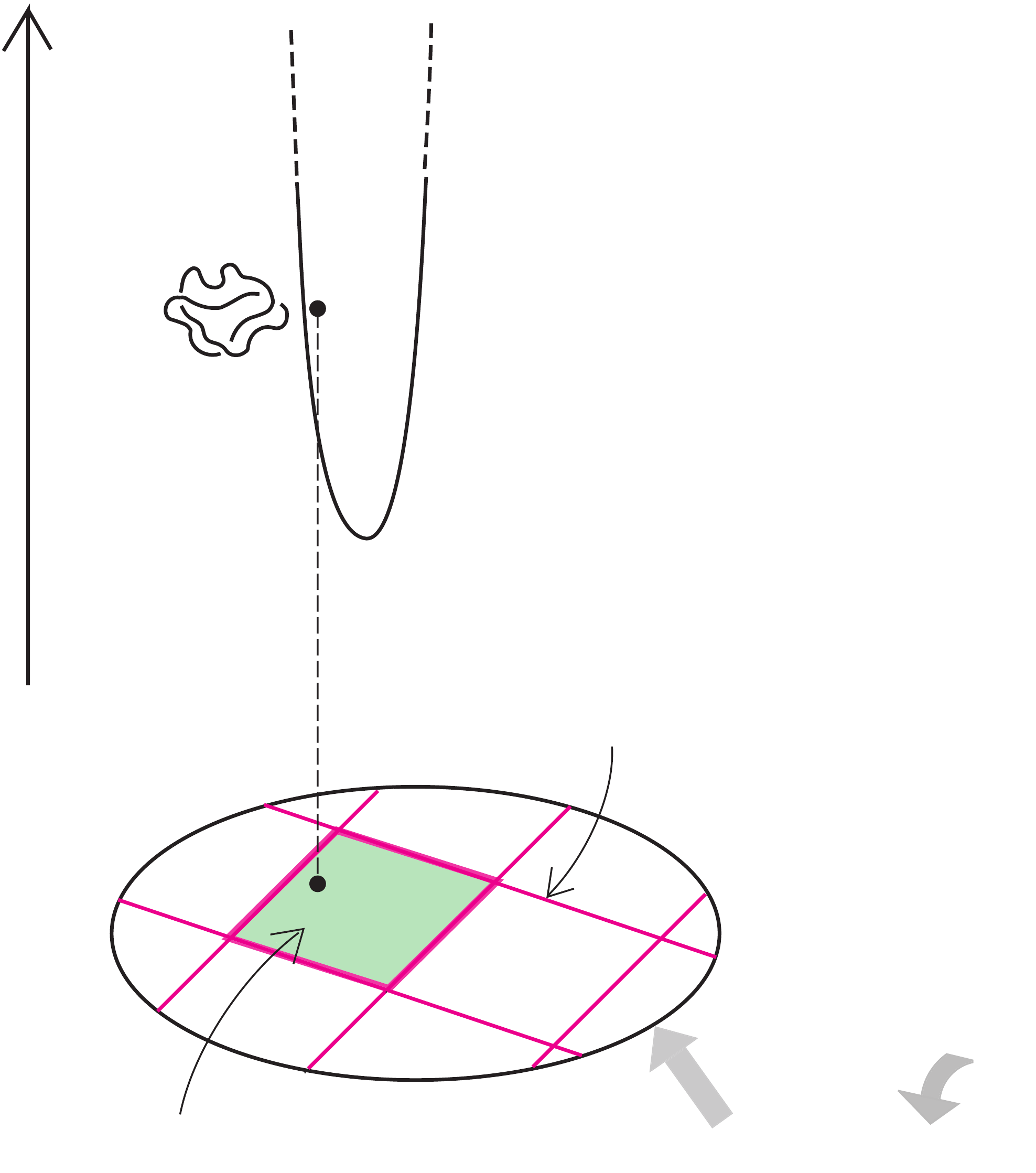}
\hskip 0.2cm
\labellist\small
\pinlabel $e$ [b] at 142 622
\pinlabel* $e([K])$ [rb] at 136 329
\pinlabel $e$--minimizer [lb] at 402 325
\pinlabel* $K_0$ [lb] at 335 120
\endlabellist
\raise20pt\hbox{\includegraphics[width=0.33\linewidth]{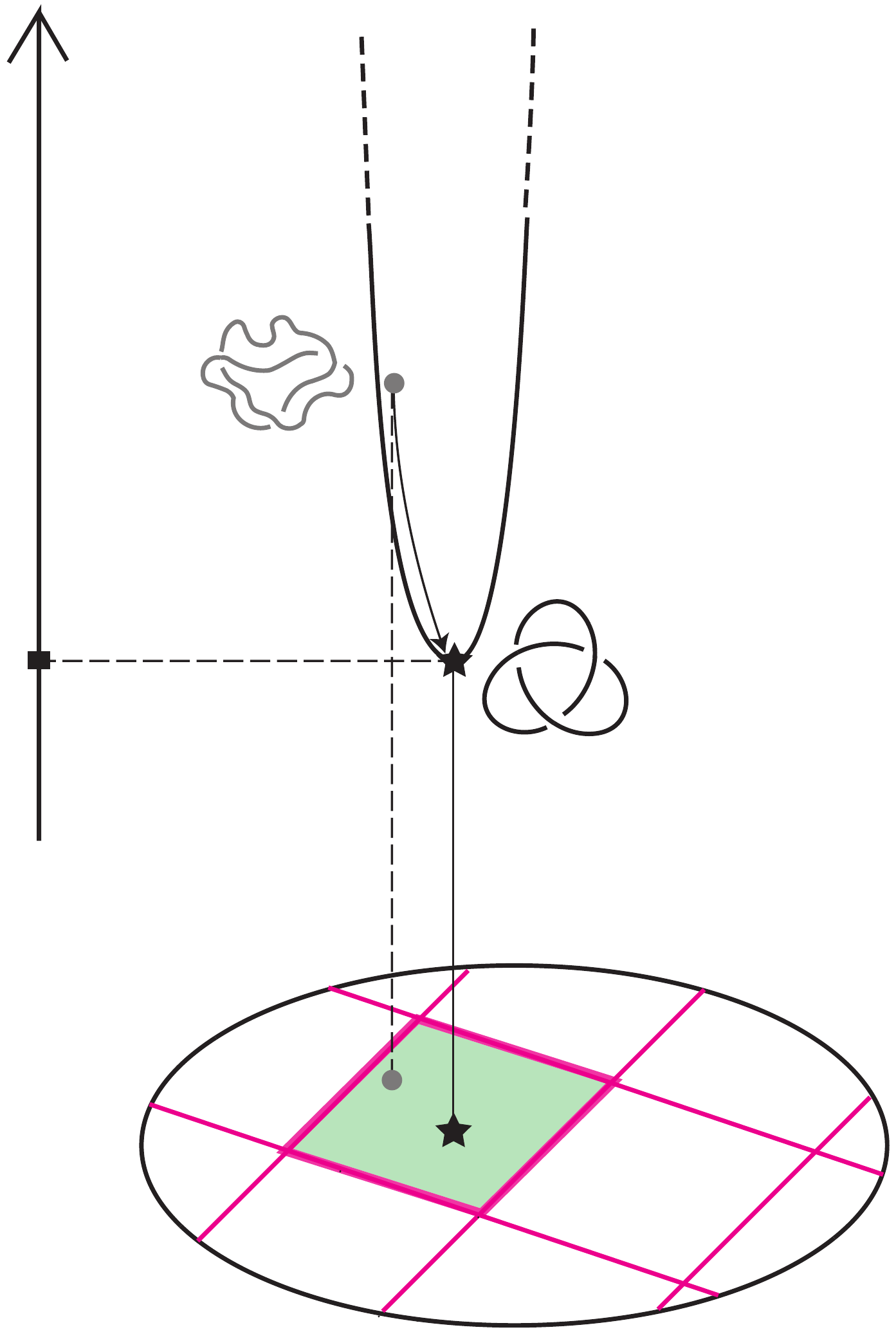}}
\end{center}
\caption{}
\label{motivation}
\end{figure}

Our strategy can be illustrated conceptually in Figures \ref{motivation}--\ref{motivation-c2ea}. 

Let $\mathcal{I}$ be the set of immersions from a circle into $\mathbb{R}^3$ (or $S^3$) and $\mathcal{D}$ the set of immersions that are not embeddings. 
Sometimes this set $\mathcal{D}$ is called the {\em discriminant} set, and an element of $\mathcal{D}$ is called a {\em singular knot}. 
Let $\mathcal{K}$ be the complement of $\mathcal{D}$ in $\mathcal{I}$, ie the space of knots. 
We will always assume that $\mathcal{I}$ is endowed with $C^2$--topology. 
Two knots $K$ and $K^{\prime}$ can be joined by a continuous path in the space of the knots $\mathcal{K}$ if and only if $K$ and $K^{\prime}$ are isotopic. 
Therefore each ``cell'' (an arcwise connected component) of $\mathcal{K}$ corresponds to an isotopy class. 

Given a knot $K$ (\fullref{motivation} (a)). 
Suppose it can be evolved along the negative gradient flow of $e$
(\fullref{motivation} (b)). 
Assume that it converges to an $e$--minimizer $K_{0}$ as time goes to infinity. 

If $K_{0}$ is isotopic to the original knot $K$ then the problem is settled. 
As the isotopy class of a knot might be changed by a crossing change, 
it should be avoided while the knot is being evolved. 
Thus we are lead to the condition below. 

\begin{defn} \rm 
We call a functional $e\co \mathcal{K}\to\mathbb{R}$ a {\em self-repulsive energy} of knots, or simply, an {\em energy of knots}, if $e(K)$ blows up as $K$ degenerates to a singular knot with double points (\fullref{motivation-c2ea}). 
\end{defn}

If $e\co \mathcal{K}\to\mathbb{R}$ is a self-repulsive energy of knots then each isotopy class is surrounded by infinitely high energy walls.

\begin{figure}[ht!]
\labellist\small
\pinlabel $e$ [b] at 18 593
\pinlabel $e$ [r] at 167 650
\pinlabel $=\infty$ [l] at 259 650
\endlabellist
\begin{center}
\includegraphics[width=0.4\linewidth]{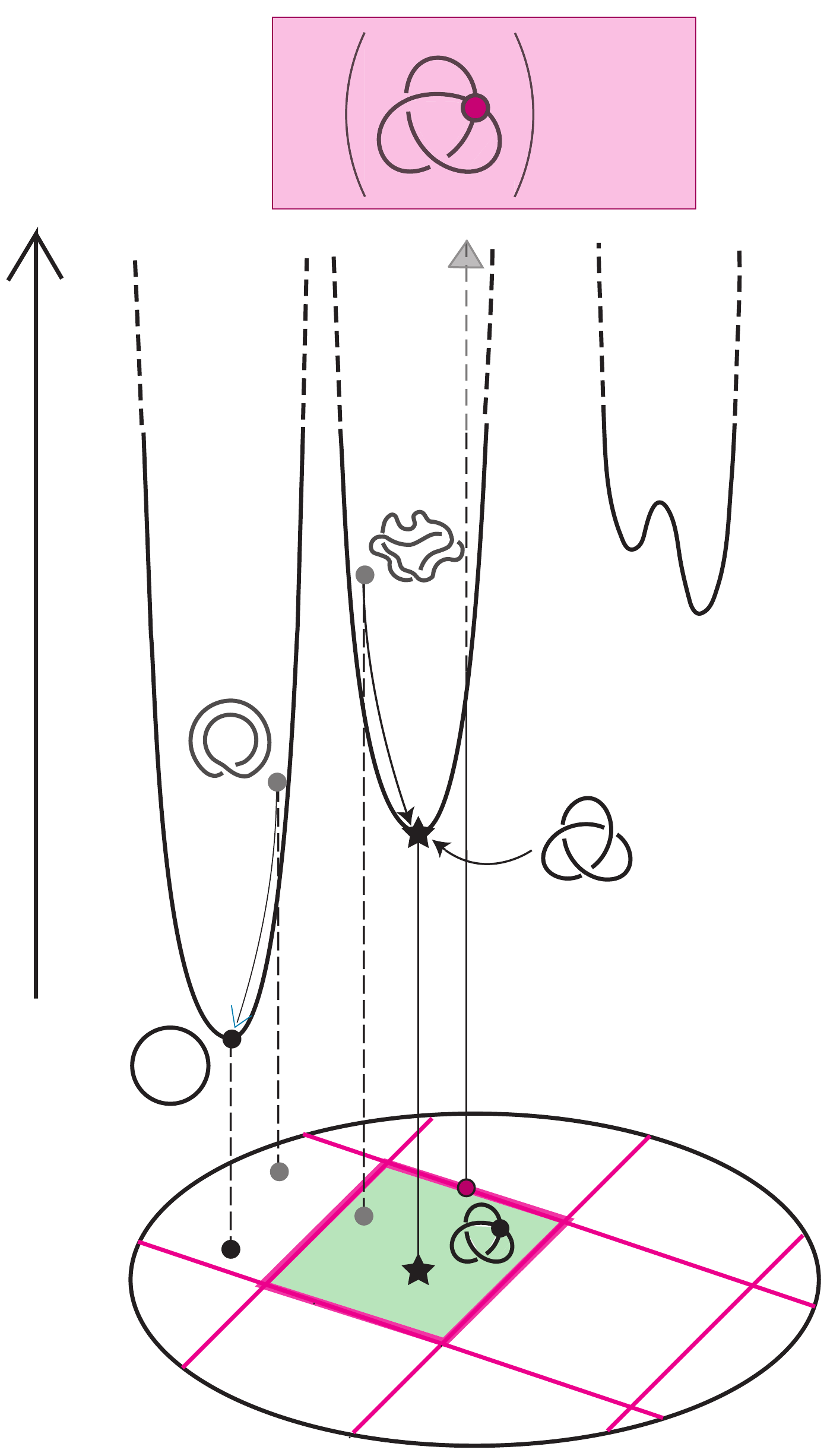}
\caption{Our energy should blow up on the discriminant set.}
\label{motivation-c2ea}
\end{center}
\end{figure}

\subsection{Renormalizations of electrostatic energy}

The first example of such an energy, $E^{(2)}_{\circ}$, was defined as the renormalization of a ``modified''  electrostatic energy of uniformly charged knots. 
The electrostatic energy of a charged knot $K$ is given by $\mbox{``$E$''}(K)=\displaystyle \int\!\!\!\int_{\!K\times K}\frac{dxdy}{|x-y|}$ which turns out to be $+\infty$ for any knot, as it blows up at the diagonal $\Delta\subset K\times K$. 
We use a trick of subtacting a function which blows up in the same order at the diagonal to produce a finite valued functional $E^{(1)}$. 
But $E^{(1)}\big(\widehat K\,\big)$ does not blow up for a singular knot $\widehat K$ with double points, ie, it is not a self-repulsive energy of knots. 
We obtain a self-repulsive energy if we make the power of $|x-y|$ in the integrand bigger than or equal to $2$. 
Let us first study the case when it is equal to $2$. 
It means that we consider the ``modified'' energy under the assumption that the magnitude of Coulomb's repulsive force  between a pair of point charges of distance $r$ is proportional to $r^{-3}$. 

\begin{defn}(O'Hara \cite{OH1})\qua
Let $d_K(x,y)$ denote the (shorter) arc-length between $x$ and $y$ (\fullref{arc_chord_length-c}). 
\begin{eqnarray}
E^{(2)}_{\circ}(K)
&=&\displaystyle{\lim_{\e \to 0}
\left\{\iint_{\!\{d_K(x,y)\ge\e\}\subset K\times K}\frac{dxdy}{|x-y|^2}
-\frac2{\,\e\,}\right\}}\label{1st_formula_of_E}\\[2mm]
&=&\displaystyle{-4+\!\iint_{\!\! K\times K}\left(\frac1{|x-y|^2}-\frac1{d_K(x,y)^2}\right)dxdy\,,}\label{2nd_formula_of_E}
\end{eqnarray}
where we assumed that the length of the knot $K$ is equal to $1$ in \ref {1st_formula_of_E}. 
\end{defn}
\begin{figure}[ht!]
\begin{center}
\begin{minipage}{.45\linewidth}
\begin{center}
\labellist\small
\pinlabel $K$ [br] at 36 214
\pinlabel $x$ [b] at 223 165
\pinlabel $y$ [t] at 160 15
\pinlabel $d_K(x,y)$ [l] at 286 127 
\pinlabel $|x-y|$ [tl] at 281 41
\endlabellist
\includegraphics[width=.6\linewidth]{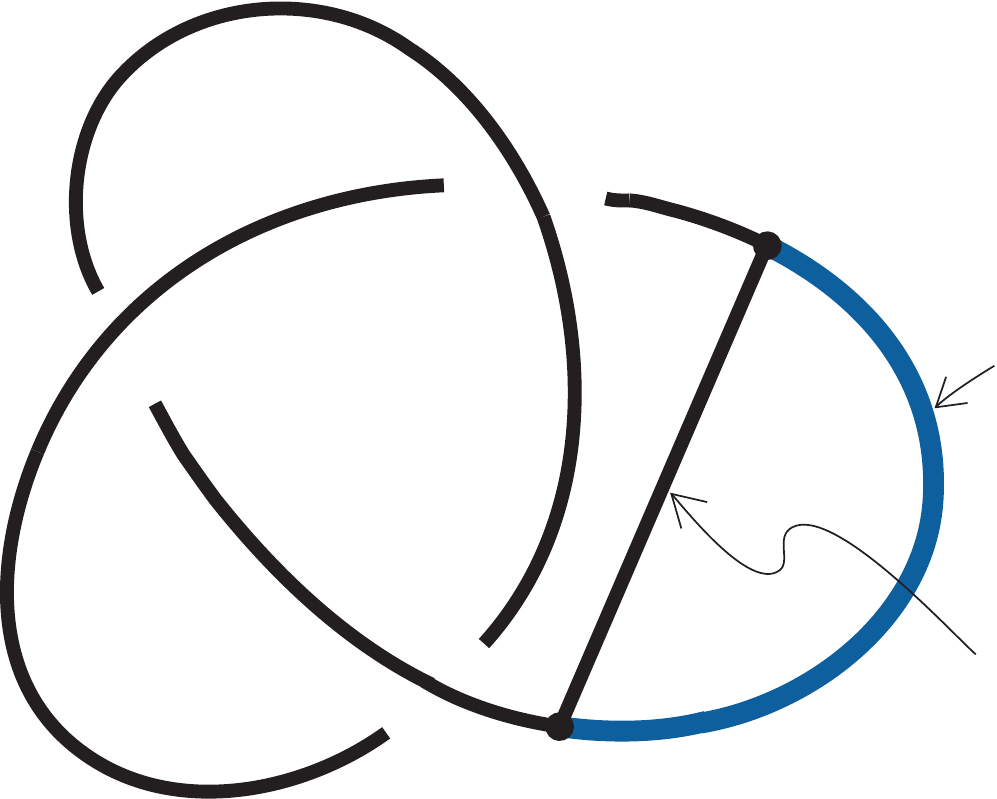}
\caption{The arc-length and the chord length}
\label{arc_chord_length-c}
\end{center}
\end{minipage}
\hskip 0.4cm
\begin{minipage}{.45\linewidth}
\begin{center}
\labellist\small
\pinlabel $K$ [br] at 27 281
\pinlabel $x$ [l] at 295 150
\pinlabel $\varepsilon$ [l] at 295 185
\pinlabel $\varepsilon$ [l] at 298 119
\pinlabel* charged [bl] at 199 300
\endlabellist
\includegraphics[width=0.5\linewidth]{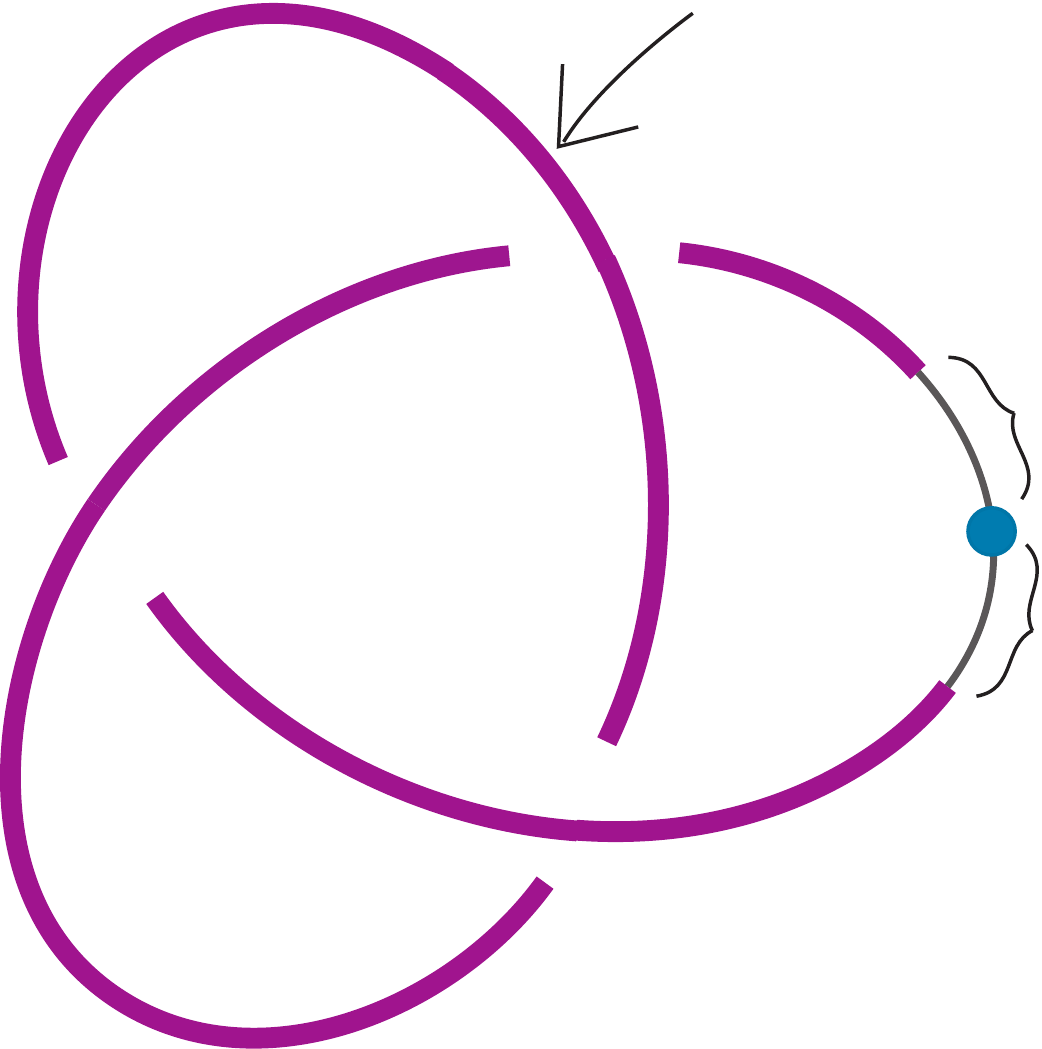}
\addtocounter{figure}1
\hypertarget{fig4}{\caption{The subarc $\{y\in K\mid d_K(x,y)\ge\e\}$}}
\label{self_avoiding_voltage-c}
\end{center}
\end{minipage}
\end{center}
\end{figure}

The term $\displaystyle \int_{\{d_K(x,y)\ge\e\}}\frac{dy}{|x-y|^2}$ in
\ref {1st_formula_of_E} expresses the ``voltage'' at point $x$ when
the subarc $\{y\in K|d_K(x,y)\ge\e\}$ is charged
(\hyperlink{fig4}{Figure 4}).  The renormalization in
\eqref{2nd_formula_of_E} can be interpreted as taking the difference
of the ``{extrinsic energy}'' based on the distance in the ambient
space (chord length) and the ``{intrinsic energy}'' based on the
distance in the knot (arc-length).  If $\Gamma_{\circ}$ denotes a
round circle then $E^{(2)}_{\circ}(\Gamma_{\circ})=0$.  It gives the
smallest value of $E^{(2)}_{\circ}$ among all knots.

If $\widetilde K$ is an open long knot, ie, the embedded line in
$\mathbb{R}^3$ which tends asymptotically to a straight line at the
both ends, its energy can be defined by dropping off the constant $-4$
(Freedman, He and Wang \cite{Fr-He-Wa}): 
$$
E^{(2)}_{\circ}(\widetilde K)=\iint_{\!\! \wtilde K\times \wtilde K}\left(\frac1{|x-y|^2}-\frac1{d_{\widetilde K}(x,y)^2}\right)dxdy. 
$$
If $\widetilde{\Gamma_{\circ}}$ is a straight line then $E^{(2)}_{\circ}(\widetilde{\Gamma_{\circ}})=0$. 

\subsection{Conformal invariance of $E^{(2)}_{\circ}$ and $E^{(2)}_{\circ}$--minimizers} 
\begin{figure}[ht!]
\begin{center}
\labellist\small
\pinlabel $C$ [t] at 100 100
\pinlabel* $r$ [r] at 116 147
\pinlabel $P$ [t] at 137 100
\pinlabel $P'$ [t] at 365 100
\pinlabel $|CP|\cdot|CP'|=r^2$ [bl] at 186 172
\endlabellist
\includegraphics[width=.38\linewidth]{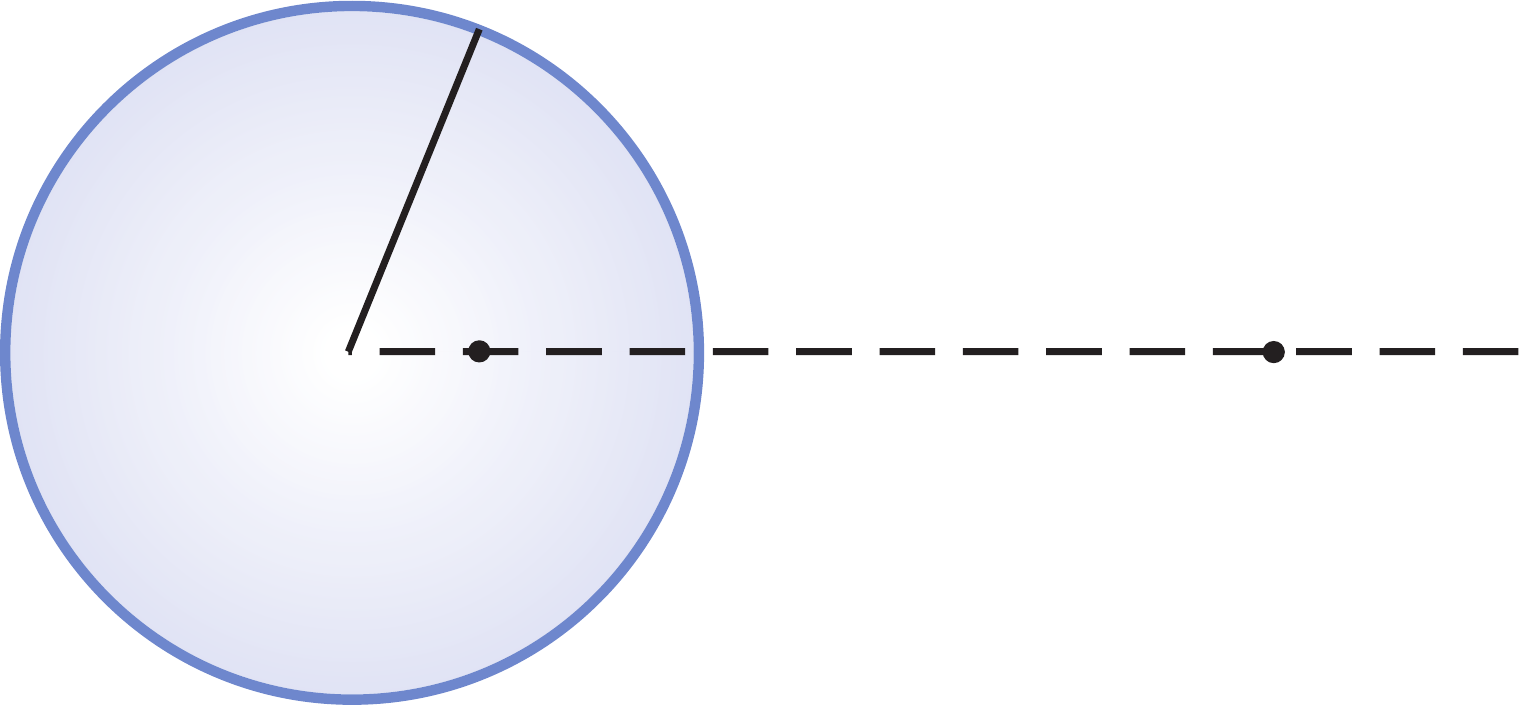}
\caption{An inversion in a sphere}
\label{inversion-c}
\end{center}
\end{figure}

The value of $E^{(2)}_{\circ}$ is invariant under rescaling or reparametrization. 
A {\em M\"obius transformation} is a transformation of $\mathbb{R}^3\cup\{\infty\}$ which can be obtained as the composition of inversions in spheres (\fullref{inversion-c}). 
\begin{thm}{\rm \cite{Fr-He-Wa}}\label{thm_fhw_mobius_inv}\qua
Let $K$ be a knot in $\mathbb{R}^3$ and $T$ a M\"obius transformation of $\mathbb{R}^3\cup\{\infty\}$. 
Then $E^{(2)}_{\circ}(T(K))=E^{(2)}_{\circ}(K)$. 
It holds even when $T(K)$ is an open long knot. 
\end{thm}
We remark that the $2$--form $\displaystyle \frac{dxdy}{|x-y|^2}$ on $K\times K\setminus\Delta$, which is an essential part of the integrand of \eqref{2nd_formula_of_E} of the definition of $E^{(2)}(K)$, is invariant under M\"obius transformations, in other words, if $K=f(S^1)$ then 
$$f^{\ast}\left(T^{\ast}\left(\frac{dxdy}{|x-y|^2}\right)\right)=f^{\ast}\left(\frac{dxdy}{|x-y|^2}\right).$$
Using the conformal invariance Freedman, He, and Wang gave a partial affirmative answer to the motivational problem: 
\begin{thm}{\rm \cite{Fr-He-Wa}}\qua
There exists an $E_{\circ}^{(2)}$--minimizer for any isotopy class of a \underline{prime} knot. 
\end{thm}

\begin{conjecture}\label{conj_Ku_Su_composite}On the other hand, Kusner and Sullivan
  \cite{Ku-Su} conjectured through numerical experiments that there
would be no $E_{\circ}^{(2)}$--minimizers in any isotopy class of a
composite knot $[K_1\sharp K_2]$, because both tangles representing
$[K_1]$ and $[K_2]$ would ``{\sl pull tight}'' to points if the knot
evolves itself to decrease its energy (\fullref{composite3}).
\begin{figure}[ht'']
\begin{center}
\includegraphics[width=.66\linewidth]{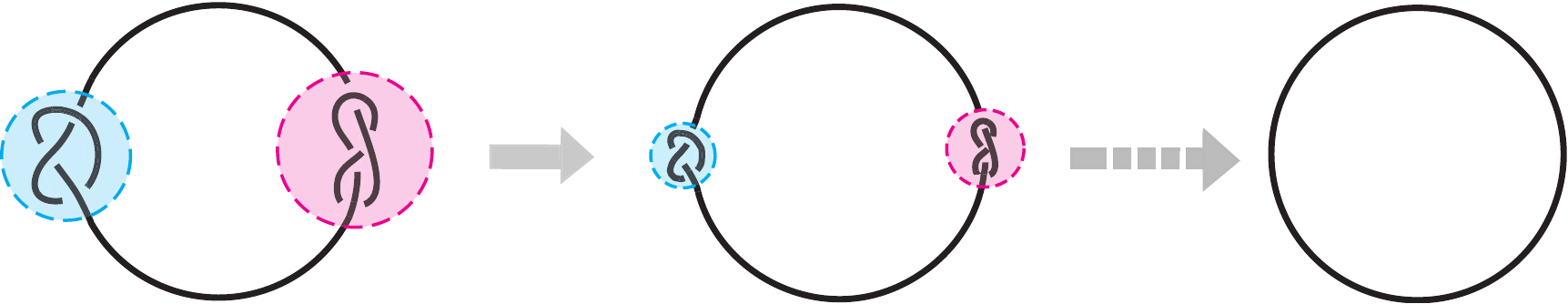}
\caption{``Pull-tight''}
\label{composite3}
\end{center}
\end{figure}

They also conjecture that 
$$E^{(2)}_{\circ}([K_1\sharp K_2])=E^{(2)}_{\circ}([K_1])+E^{(2)}_{\circ}([K_2]).$$ 
\end{conjecture}

\begin{figure}[ht!]
\begin{center}
\includegraphics[width=.66\linewidth]{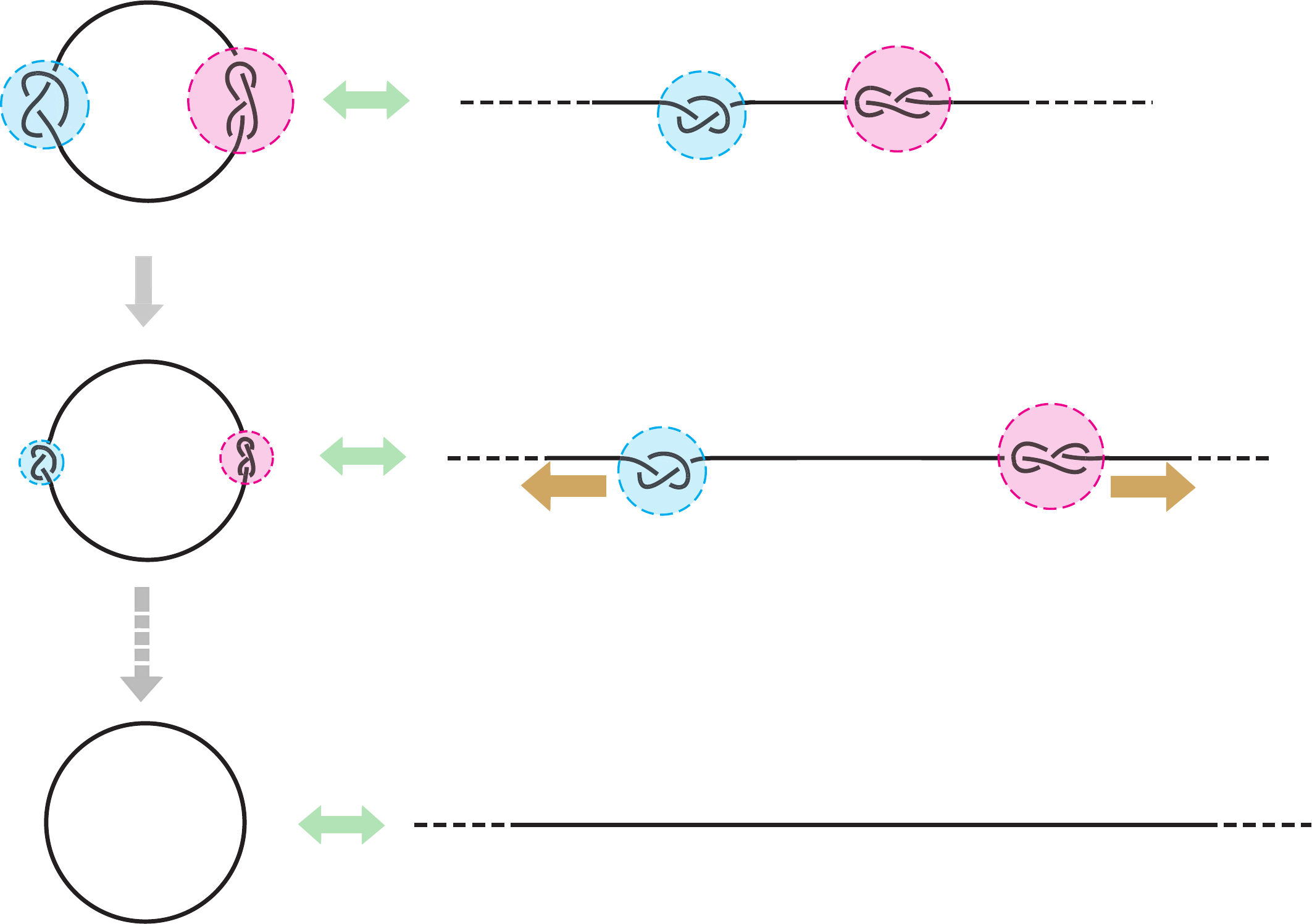}
\caption{``Pull-tight'' in open long knots}
\label{composite}
\end{center}
\end{figure}

This can be explained as follows. 
Consider an open long knot by an inversion in a sphere with center on the knot. 
If a knot pulls tight then the two tangles in the open long knot corresponding to $[K_1]$ and $[K_2]$ move in the opposite ways so that they become more and more distant from each other (\fullref{composite}).  
As the distance between the two tangles tends to $+\infty$ the interaction between them in the integral of $E^{(2)}_{\circ}$ tends to $0$. 

\begin{figure}[ht!]
\begin{center}
\labellist\small
\pinlabel $K_1$ [bl] at 131 377
\pinlabel $K_2$ [bl] at 337 377
\pinlabel $K_{\infty}$ [bl] at 585 377
\pinlabel {If the knot is prime} at 176 217
\pinlabel {$E^{(2)}_{\circ}$ is continuous} <0pt, 7pt> at 523 202
\pinlabel {wrt $C^2$--topology} <0pt, -7pt> at 523 202
\endlabellist
\includegraphics[width=.55\linewidth]{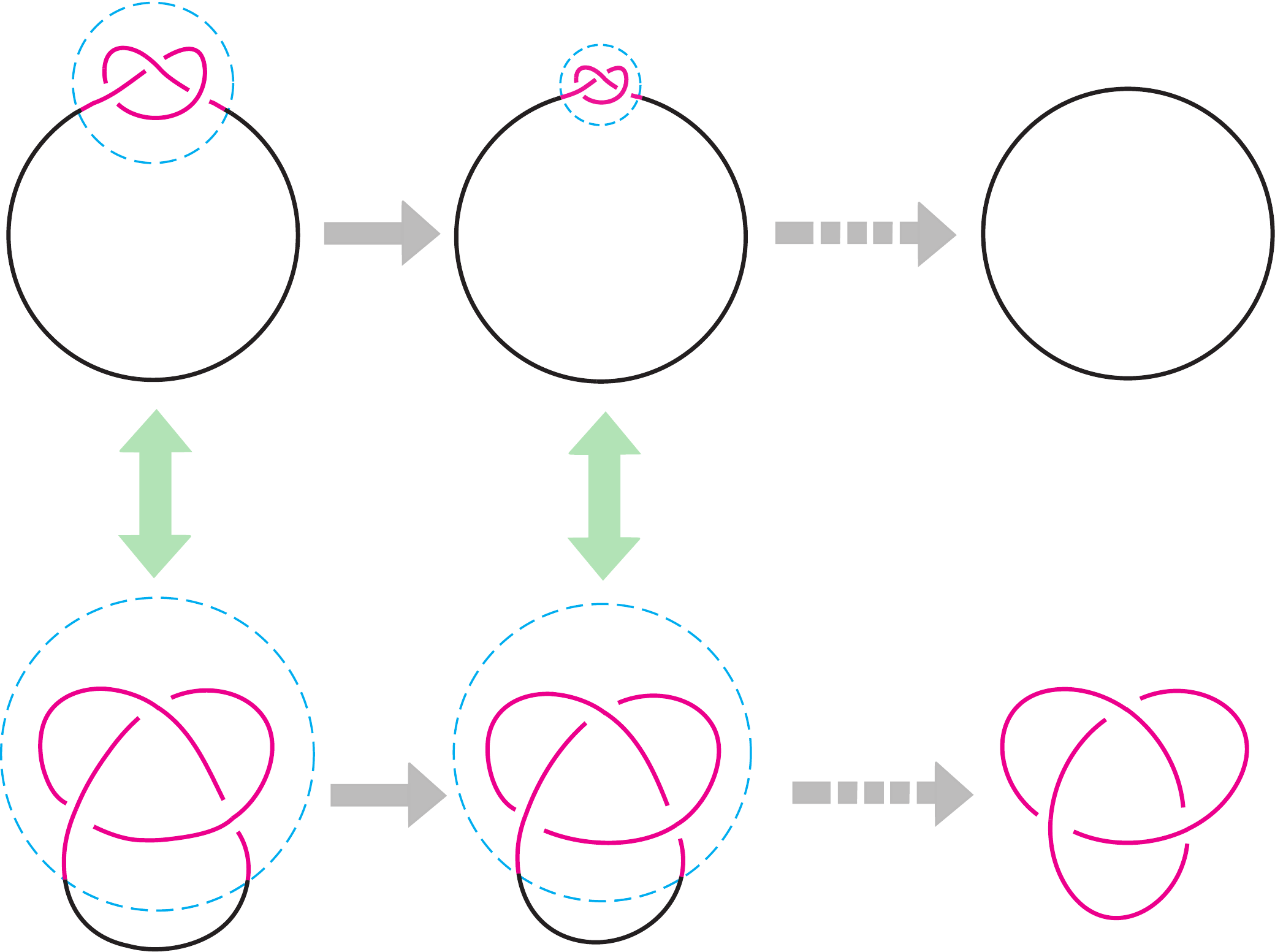}
\caption{In the case of an isotopy class of a prime knot}
\label{pull-tight-mobius}
\end{center}
\end{figure}

In the case of prime knots, the pull-tight can be avoided as follows. 
Suppose $\{K_n\}$ is a sequence of knots in an isotopy class of a prime knot $[K]$ with $\displaystyle \lim_{n\to\infty}E^{(2)}_{\circ}(K_n)=E^{(2)}_{\circ}([K])$. 
Applying M\"obius transformations if necessary, we can obtain a new sequence of ``{\sl relaxed}\,'' knots $\{K^{\p}_n\}\subset [K]$ so that the pull-tight does not occur and hence $\displaystyle \lim_{n\to\infty}K^{\p}_n$ belongs to the same isotopy class $[K]$. 
As $E^{(2)}_{\circ}(K^{\p}_n)=E^{(2)}_{\circ}(K_n)$ the limit $\displaystyle \lim_{n\to\infty}K^{\p}_n$ is an $E^{(2)}_{\circ}$--minimizer of $[K]$ (\fullref{pull-tight-mobius}). 

Here are some remarks: 

\begin{enumerate}
\item Z-X He \cite{He} showed that $E_{\circ}^{(2)}$--minimizers are
smooth.
\item Suppose $K$ is an $E_{\circ}^{(2)}$--minimizer of an isotopy class $[K]$. 
Then, for any M\"obius transformation $T$, at least one of $T(K)$ and its mirror image $T(K)^{\ast}$ belongs to $[K]$. 
Since $E_{\circ}^{(2)}(T(K))=E_{\circ}^{(2)}(T(K)^{\ast})=E_{\circ}^{(2)}(K)$ it follows that $T(K)$ or $T(K)^{\ast}$ is an $E_{\circ}^{(2)}$--minimizer of $[K]$. 
Therefore, there are uncountably many $E_{\circ}^{(2)}$--minimizers for each isotopy class of a non-trivial prime knot. 
\item The (cardinal) number of $E_{\circ}^{(2)}$--minimizers of an isotopy class of a prime knot modulo the action of the M\"obius group is not known. 
\item It is not known whether there exists an $E_{\circ}^{(2)}$--critical unknot which is not a round circle. 
If not, it implies Hatcher's results \cite{Ha} that the set of unknots in $S^3$ deformation retracts onto the set of great circles. 

Numerical experiments show that $E_{\circ}^{(2)}$ can untie Ochiai's
unknot (Kauffman, Huang and Greszczuk \cite{Ka-Hu-Gr}) and
``Freedman's unknot'' (Kusner and Sullivan \cite{Ku-Su}); see 
\fullref{Freedman_unknot}.
   \begin{figure}[ht!]
   \cl{\includegraphics[width=.25\linewidth]{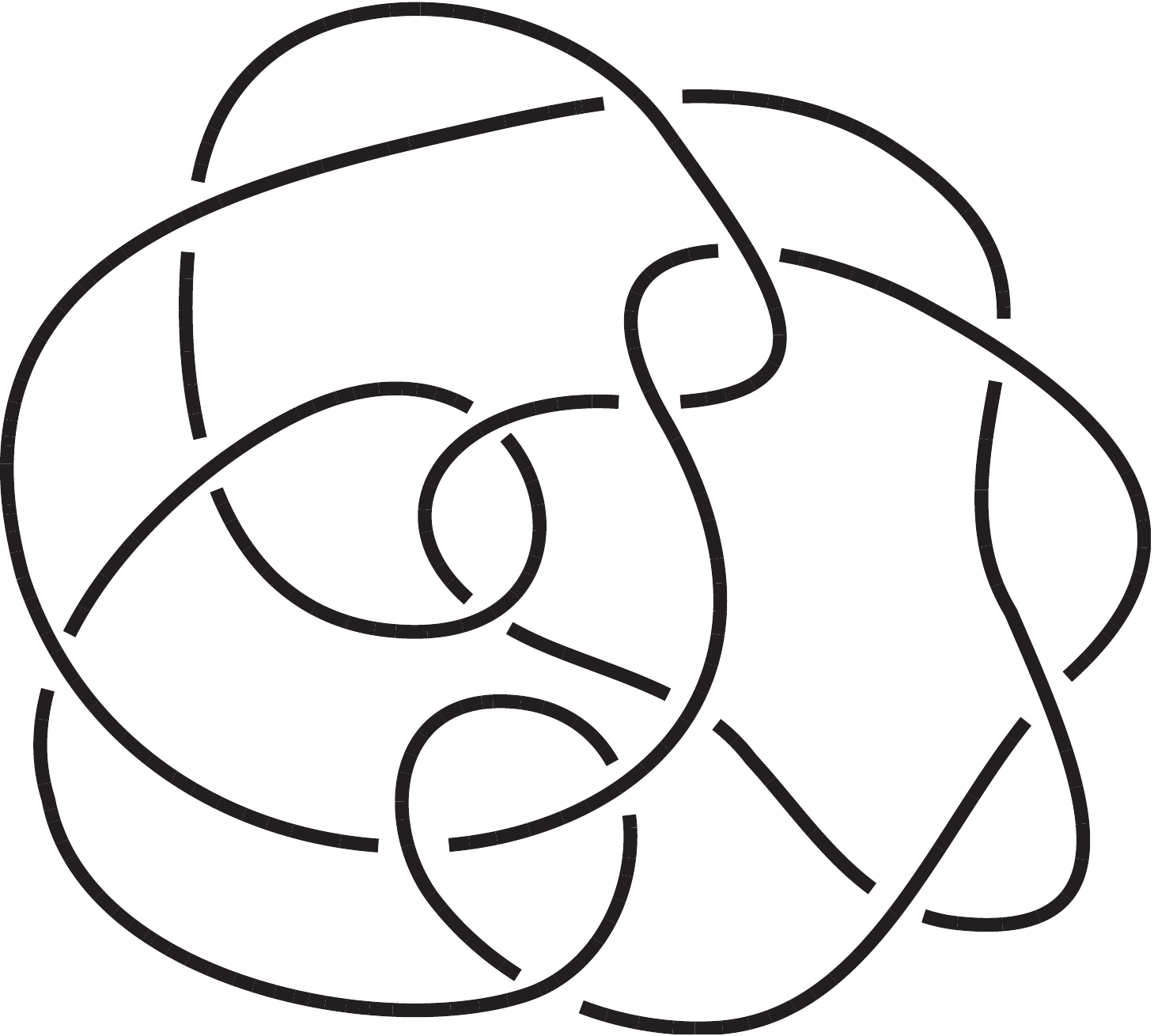}\hspace{1cm}
   \includegraphics[width=.5\linewidth]{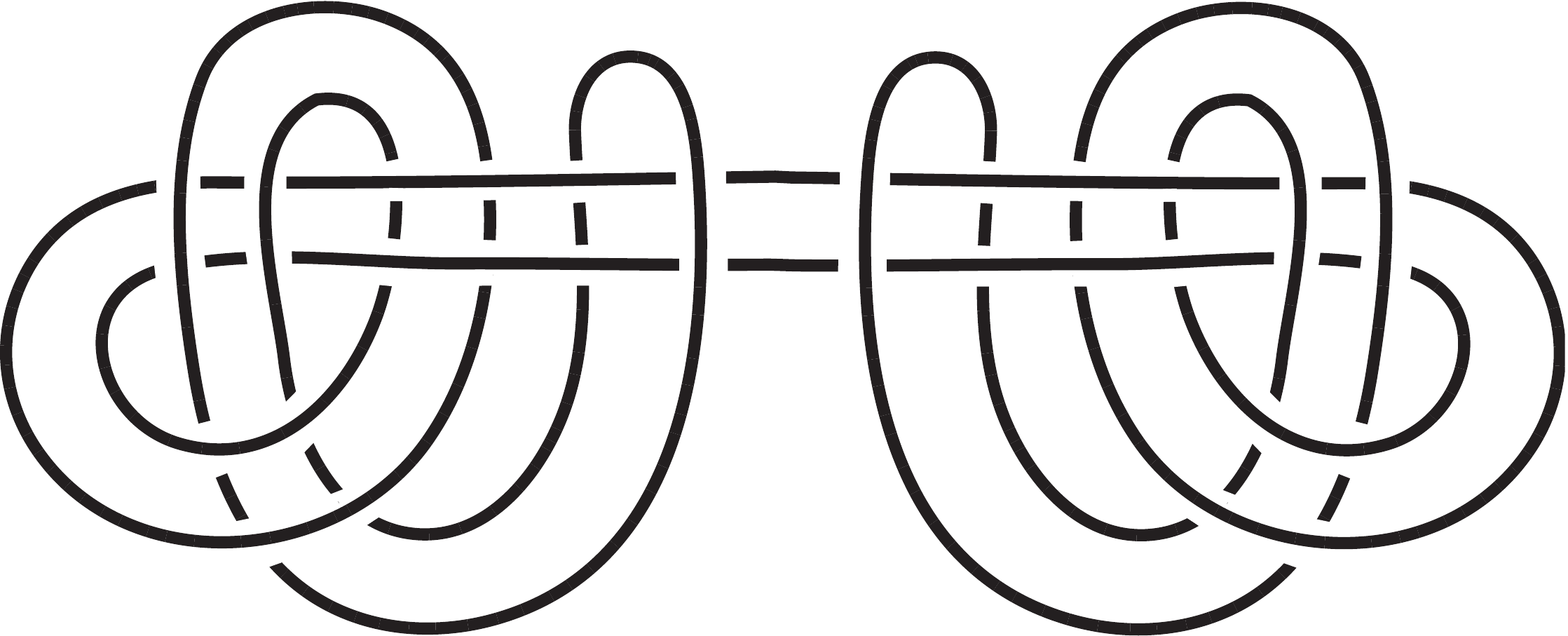}}
   \cl{\small\hglue -1cm (a)\hspace{5.5cm}(b)\hspace{.5cm}}   
\caption{(a)\qua Ochiai's unknot \quad (b)\qua Freedman's unknot}
   \label{Freedman_unknot}
   \end{figure}
\item Using numerical experiments, Kusner and Sullivan conjecture that there exist unstable critical points in the isotopy class of a $(p,q)$ torus knot if both $p$ and $q$ are greater than $2$. 
\item There are no known minimum values of $E_{\circ}^{(2)}$ of an isotopy class of a non-trivial knot which are obtained theoretically, like $6\pi^2$. 
\item It is an open problem whether $E^{(2)}_{\circ}$--minimizers are isolated in $\mbox{Emb}(S^1,\mathbb{R}^3)/\sim$, where $\sim$ is generated by M\"obius transformations and reparametrizations. 
\end{enumerate}

Various kinds of generalization of $E_{\circ}^{(2)}$ have been studied. 
\begin{enumerate}
\item $E_{\circ}^{(2)}$ can be defined for a link $L=K_1\cup\cdots\cup K_n$ \cite{Fr-He-Wa}. 
We do not need renormalization for the cross term $\displaystyle E(K_i, K_j)=\int_{x\in K_i}\int_{y\in K_j}\frac{dxdy}{|x-y|^2}$ $(i\ne j)$. 
\item Similar energies are studied by Buck, Orloff and Simon \cite{Bu-Or,Bu-Si}. 
The integrands are the products of the $2$--form $\frac{dxdy}{|x-y|^2}$ and functions which kill the explosion of the integral at the diagonal. 
\item A conformally invariant energy for surfaces was studied by
  Auckly and Sadun \cite{Au-Sa}. 
\item A conformally invariant energy for hypersurfaces using
  conformally defined angles was studied in Kusner and Sullivan \cite{Ku-Su}. 
\item Fixing a knot $K$, $\displaystyle \iint_{\!\!K\times
  K}|x-y|^sdxdy$ can be considered a complex valued function of a
  complex variable $s$ (Brylinski \cite{Bry}). 
\end{enumerate}

\subsection{Generalization to produce energy minimizers}
\fullref{conj_Ku_Su_composite} implies that $E_{\circ}^{(2)}$ does not give a completely affirmative solution to our motivational \fullref{motivation_probl}. 
We have two ways to generalize $E_{\circ}^{(2)}$ so that all the isotopy classes have energy minimizers. 

One is to make the power of $|x-y|$ in the integrand bigger than $2$, and the other is to change the metric of the ambient space. 
In each case, our energies are no longer conformally invariant. 

\begin{defn} \rm 
Let $K$ be a knot with total length $1$. 
Put 
$$E^{(\a)}(K)=\!\int\!\!\!\!\int_{\! K\times K}\!\left(
\frac{\,1\,}{|x-y|^{\a}}-\frac{\,1\,}{d_K (x,y)^{\a}}\right) \!dxdy. $$
\end{defn}
\begin{thm}{\rm (O'Hara \cite{OH2,OH3})}\qua
$E^{(\a)}$ is well-defined if $\a<3$, and is self-repulsive if $\a\ge2$. 
There exists an $E^{(\a)}$--minimizer for any isotopy class if $\a>2$. 
\end{thm}

Let $M$ be a Riemannian manifold. 
Define 
$$\begin{array}{l}
d_{M}(x,y)=\inf\{\mbox{Length of path joining $x$ and $y$}\},\\[1mm]
E^{(\a)}_{M}(K)\!=\!\!\displaystyle \int\!\!\!\int_{K\times K}\!\left(\frac{\,1\,}{d_M(x,y)^{\a}}-\frac{\,1\,}{d_K(x,y)^{\a}}\right)\!dxdy. 
\end{array}
$$

\begin{thm}{\rm (O'Hara \cite{OH4})}\qua
Let $M$ be a compact manifold.  Then there exists an
$E^{(\a)}_{M}$--minimizer for any isotopy class if $\a>2$.
\end{thm}

We conjecture that the Theorem above also holds for $\a=2$ if $M=S^3$. 

\subsection{Related topics}
Energy of knots gave rise to {\sl geometric knot theory}, in which we study functionals to measure how complicated a knot is embedded and look for ``optimal knots'' with respect to those functionals. 

One of the functionals which are intensively studied recently is the {\em rope length} (Cantarella, Kusner, Sullivan, Stasiak, et al.), which measures how long a rope of unit diameter is needed to make a given knot, or its equivalents, thickness (Buck, Rawdon, Simon, et al.) and global radius of curvature (Gonzalez, Maddocks, Smutny). 

\section{A viewpoint from conformal geometry}

{\it This is joint work with R\'emi Langevin.}

\medskip
We can give a new interpretation of $E_{\circ}^{(2)}$ using what is invariant under M\"obius transformations, such as circles, spheres, and angles. 

\subsection{Minkowski space}
The {\em Minkowski space} $\mathbb{R}^5_1$ is $\mathbb{R}^5$ with the non-degenerate indefinite quadratic form with index $1$: 
$$\langle \vect x, \vect x \rangle\!=\!-x_0{}^2+x_1{}^2+\cdots +x_4{}^2.$$
The set of linear isomorphisms which preserve the Lorentz metric is called the {\em Lorentz group}: 
$$
O(4,1)=\left\{A\in GL(5,\mathbb{R})\,\left|\,{}^t\!AJ^{\,5}_{1}A=J^{\,5}_{1}\right\}\right.\!,
\mbox{ where }J^{\,5}_{1}=
\left(
\begin{array}{cccc}
-1 &   &   &    \\
 &  1 &   &  \!\!{\smash{\raise10pt\hbox{\Large  O}}}\!\\
 &   &  \ddots &   \\
\!{\smash{\hbox{\Large  O}}}\!\! &     &   & 1 
\end{array}\right).
$$ 
A non-zero vector $\vect v$ in $\mathbb{R}^{5}_{1}$ is called {\em
spacelike} if $\langle\vect v,\vect v\rangle>0$, {\em lightlike} if
$\langle\vect v,\vect v\rangle=0$ and $\vect v\ne\vect 0$, and {\em
timelike} if $\langle\vect v,\vect v\rangle<0$.  The set of lightlike
vectors and the origin $V=\left.\left\{\vect v\in\mathbb{R}^{5}_1 \,
\right|\, \langle\vect v,\vect v\rangle=0\right\}$ is called the {\em
light cone}.
The hyperquadric $\varLambda=\{\vect v\in\mathbb{R}^{5}_1|\langle\vect v,\vect v\rangle=1\}$ is called the {\em de Sitter space}. 
\begin{figure}[ht!]
\begin{center}
\labellist\small
\pinlabel $0$ [br] <-5pt,0pt> at 289 298
\pinlabel* $x_1$ [tr] at 76 209
\pinlabel $x_4$ [l] at 479 298
\pinlabel $x_0$ [b] at 289 628
\pinlabel* {hyperbolic space $\mathbb{H}^3$} [bl] at 402 562
\pinlabel $\mathbb{R}^5_1$ at  42 601
\pinlabel* $S^3(1)$ [r] <0pt,-2pt> at 136 394
\pinlabel {de Sitter space $\varLambda$} [l] at 404 227
\pinlabel {light cone $V$} [l] at 487 175
\endlabellist
\includegraphics[width=.5\linewidth]{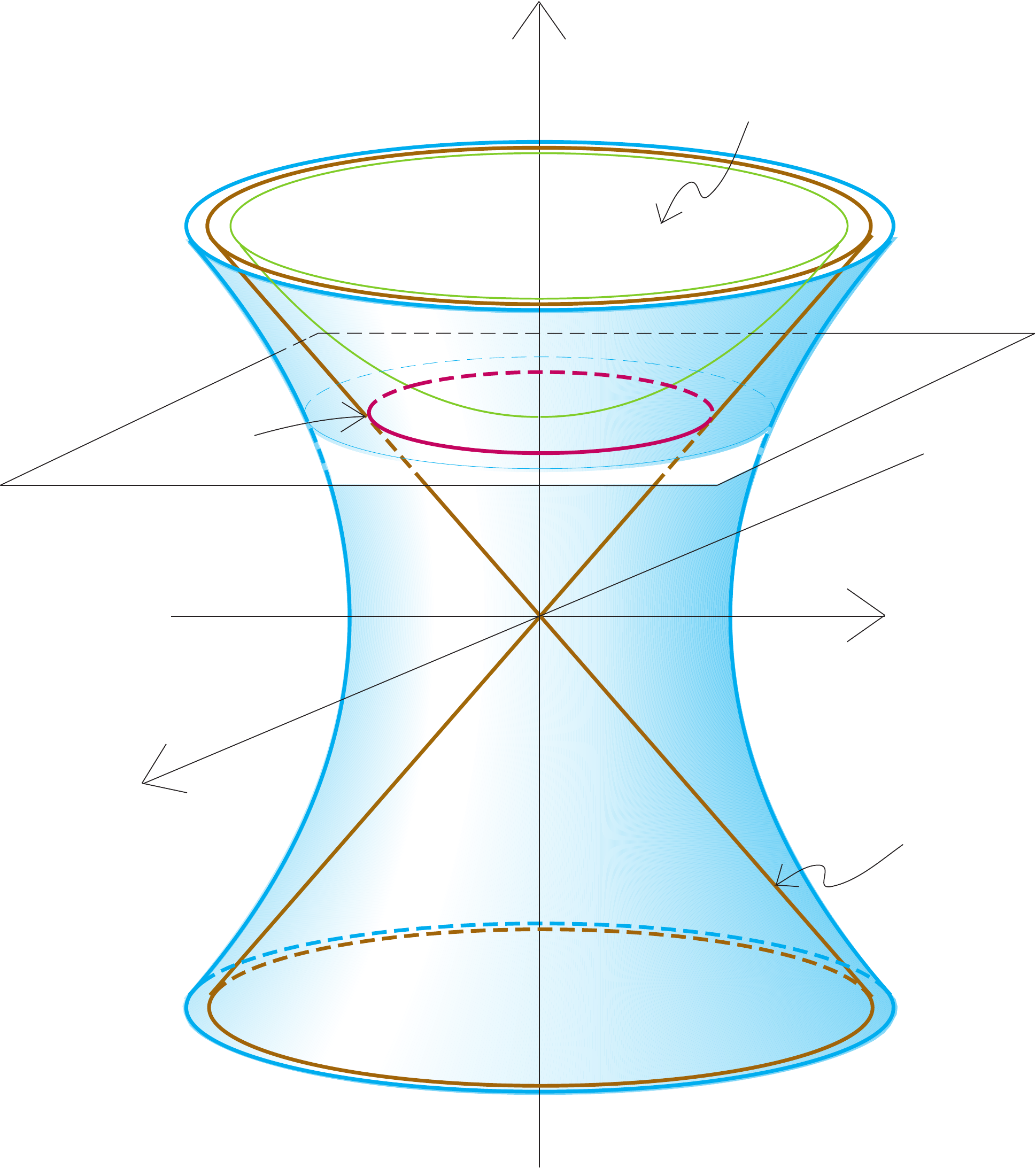}
\caption{}
\label{quadric-ce3}
\end{center}
\end{figure}
The $3$--sphere $S^3$ can be realized in $\mathbb{R}^{5}_1$ as the set of lines through the origin in the light cone $V=\left\{\langle\vect v,\vect v\rangle=0\right\}$. 
We will denote it by $S^3(\infty)$. 
It can also be identified with the intersection of the light cone and the hyperplane $\{\vect x\in\mathbb{R}^{5}_1 \,|\, x_{0}=1\}$: 
$$
S^{3}(1)=\left\{(1, x_1, x_2, x_3, x_{4})\, \left|\, x_1{}^2+x_2{}^2+x_3{}^2+x_{4}{}^2=1\right\}.\right.
$$
The Lorentz group $O(4,1)$ acts on $V$ and $\varLambda$. 
It also acts transitively on $S^3$ as the action on the set of lines in the light cone. 
This action is called a {\em M\"obius transformation}. 

\subsection{de Sitter space as the set of spheres}\label{sec_de_Sitter_sp_as_S}
Put $\mathcal{S}(2,3)=\{\varSigma \,|\, \mbox{an oriented $2$--sphere in }S^3\}$. 
Then there is a bijection between $\mathcal{S}(2,3)$ and the de Sitter space $\varLambda$. 
\begin{figure}[ht!]
\begin{center}
\labellist\small
\pinlabel $\varSigma$ [b] at 213 454
\pinlabel $\varPi$ [bl] at 392 544
\pinlabel $\sigma$ [t] at 255 283
\pinlabel $l=\varPi^{\perp}$ [bl] at 376 348
\pinlabel $\varLambda$ [b] at 391 100
\pinlabel {light cone} [l] at 364 179
\endlabellist
\includegraphics[width=.4\linewidth]{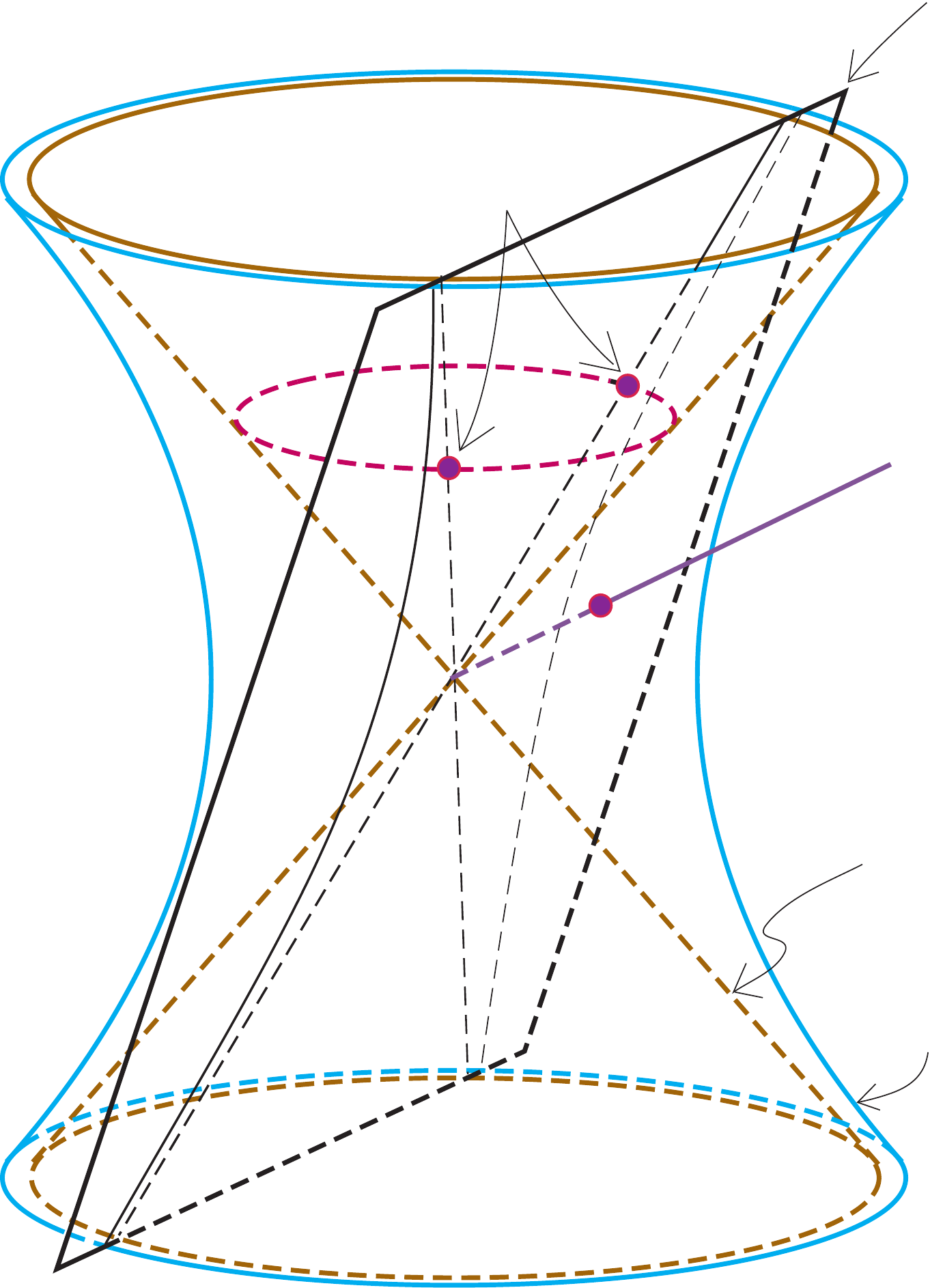}
\caption{The bijection between $\mathcal{S}(2,3)$ and $\varLambda$}
\label{Lambda-c-e}
\end{center}
\end{figure}
Let $\varSigma$ be an oriented $2$--sphere in $S^3$. 
In the Minkowski space $\mathbb{R}^5_1$, $\varSigma$ can be realized as the intersecton of $S^3$ and an oriented $4$--dimensional subsapce $\varPi$ through the origin (\fullref{Lambda-c-e}). 
Let $\sigma\in\varLambda$ be the endpoint of the positive unit normal vector to $\varPi$. 
Then the map $\varphi\co  \mathcal{S}(2,3)\ni\varSigma\mapsto \sigma\in\varLambda$ is the bijection we want. 
Moreover, since this bijection is defined only by means of the pseudoorthogonality, it is preserved under the action of the Lorentz group $O(4,1)$, ie, $\varphi(A\cdot\varSigma)=A\varphi(\varSigma)$ for $A\in O(4,1)$. 

\subsection{Willmore Conjecture}
In his attempt to solve the Willmore Conjecture (stated below)
Langevin has been interested in defining conformally invariant
functionals on the space of surfaces and knots by means of integral
geometry in the Minkowski space.  One of his functionals turned out to
be a self-repulsive energy of knots (Langevin and O'Hara \cite{LO1}).
That was the beginning of our joint work.  Let us make a short comment
on the Willmore Conjecture.

Let $\iota\co T^2\to \mathbb{R}^3$ be a smooth embedding. 
Let $\kappa_1,\kappa_2$ be the principal curvatures. 
The Willmore functional $W$ is defined by 
$$W(\iota)\!=\!\!\int_{T^2}\!\left(\frac{\kappa_1+\kappa_2}{2}\right)^2\!dv
\!=\!\!\int_{T^2}\!\left(\frac{\kappa_1-\kappa_2}{2}\right)^2\!dv. $$
The second equality comes from the Gauss-Bonnet theorem. 
The integrand of the right hand side is known to be invariant under M\"obius transformations, so is $W$. 
\begin{wconjecture}\rm 
We have $W(\iota)\ge 2\pi^2$. 
The equality holds if and only if $\iota(T^2)$ is a torus of revolution $T_{\sqrt2,1}$ modulo M\"obius transformations, where $T_{\sqrt2,1}$ can be obtained by rotating around the $z$--axis a circle with radius $1$ in the $xz$--plane whose center is distant form the $z$--axis by $\sqrt2$. 
\end{wconjecture}
One of the important contributions is the following theorem due to Bryant. 
Fix a stereographic projection $\pi\co S^3\to\mathbb{R}^3\cup\{\infty\}$. 
\begin{thm}{\rm (Bryant \cite{Br})}\qua
Let $\varSigma_{\frac{2}{\kappa_1(p)+\kappa_2(p)}}(p)$ be a sphere which is tangent to $\iota(T^2)$ at $p$ with curvature $\frac{\kappa_1(p)+\kappa_2(p)}{2}$. 
Define $\psi_{\iota}\co  T^2\to\varLambda$ by 
$$\psi_{\iota}(p)= \varphi\circ\pi^{-1}\left(\varSigma_{\frac{2}{\kappa_1(p)+\kappa_2(p)}}(p)\right),$$ 
where $\varphi\co \mathcal{S}(2,3)\to\varLambda$ is the bijection given in {\rm \fullref{sec_de_Sitter_sp_as_S}}. 
Then $W(\iota)$ is equal to the area of $\psi_{\iota}(T^2)$. 
\end{thm}
The area is given with respect to the indefinite metric of $\varLambda\subset\mathbb{R}^5_1$, which is $SO(4,1)$--invariant. 
It does not depend on the stereographic projection $\pi$. 

\subsection{Infinitesimal cross ratio} 
Let us introduce the {\sl infinitesimal cross ratio} $\Omega$, which plays an important role in our study. 
It is a complex valued $2$--form on $K\times K\setminus \Delta$. 
(The imaginary part might not be smooth.) 
It is conformally invariant. 

We explain that its real part can be interpreted in two ways. 
They correspond to two kinds of interpretations of $S^3\times S^3\setminus\Delta$ which contains $K\times K\setminus\Delta$. 
One is as the total space of the cotangent bundle $T^{\ast}S^3$, which enables us to consider $\R \Omega_K$ as the pull-back of the canonical symplectic form of the cotangent bundle $T^{\ast}S^3$. 
The other is as the set of oriented $0$--spheres in $S^3$ which has a natural pseudo-Riemannian structure coming from that of the Minkowski space, which enables us to consider $\R \Omega_K$ as a signed area form. 

Let us begin with a geometric definition of the infinitesimal cross ratio. 

\begin{figure}[ht!]
\cl{\hspace{-0.7cm}
\labellist\small
\pinlabel $K$ [b] at 393 308
\pinlabel $x$ [b] at 238 179
\pinlabel $x{+}dx$ [l] at 198 151
\pinlabel* $y$ [tl] at 288 103
\pinlabel* $y{+}dy$ [l] at 351 147
\pinlabel $\varSigma(x,x{+}dx,y,y{+}dy)$ [b] at 179 293
\endlabellist
\includegraphics[height=4.2cm]{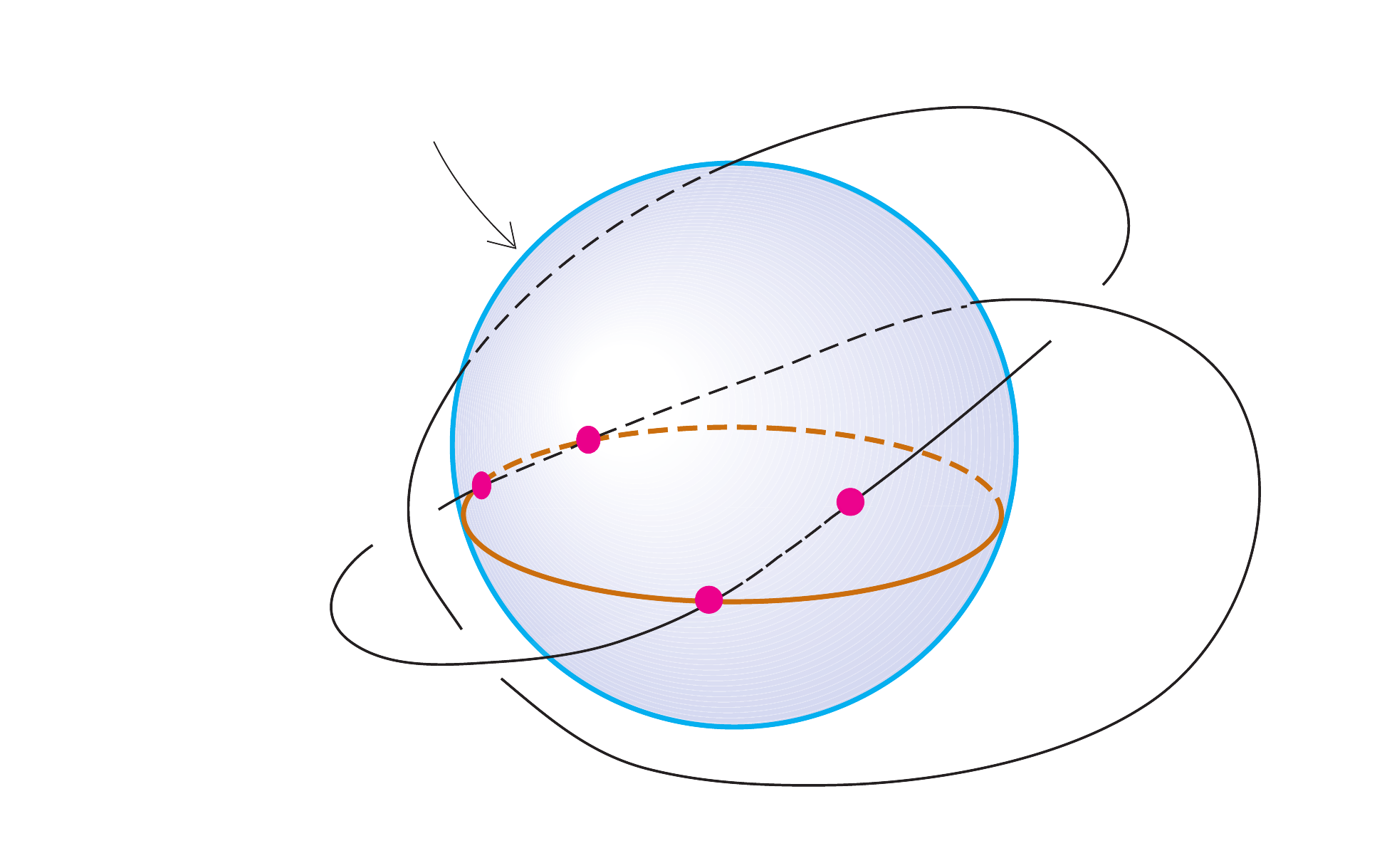}
\labellist\small
\pinlabel $K$ [b] at 393 308
\pinlabel $x$ [b] at 238 179
\pinlabel $x{+}dx$ [l] at 198 151
\pinlabel* $y$ [tl] at 288 103
\pinlabel* $y{+}dy$ [l] at 351 147
\pinlabel $\varSigma(x,x{+}dx,y,y{+}dy)$ [b] at 179 293
\pinlabel $\varPi\cong\mathbb{C}$ [rb] at 92 139
\pinlabel $C(x,x,y)$ [t] <10pt, 0pt> at 449 31
\pinlabel $\tilde y{+}\wwtilde{dy}$ [l] <4pt, 0pt> at 431 108
\endlabellist
\includegraphics[height=4.2cm]{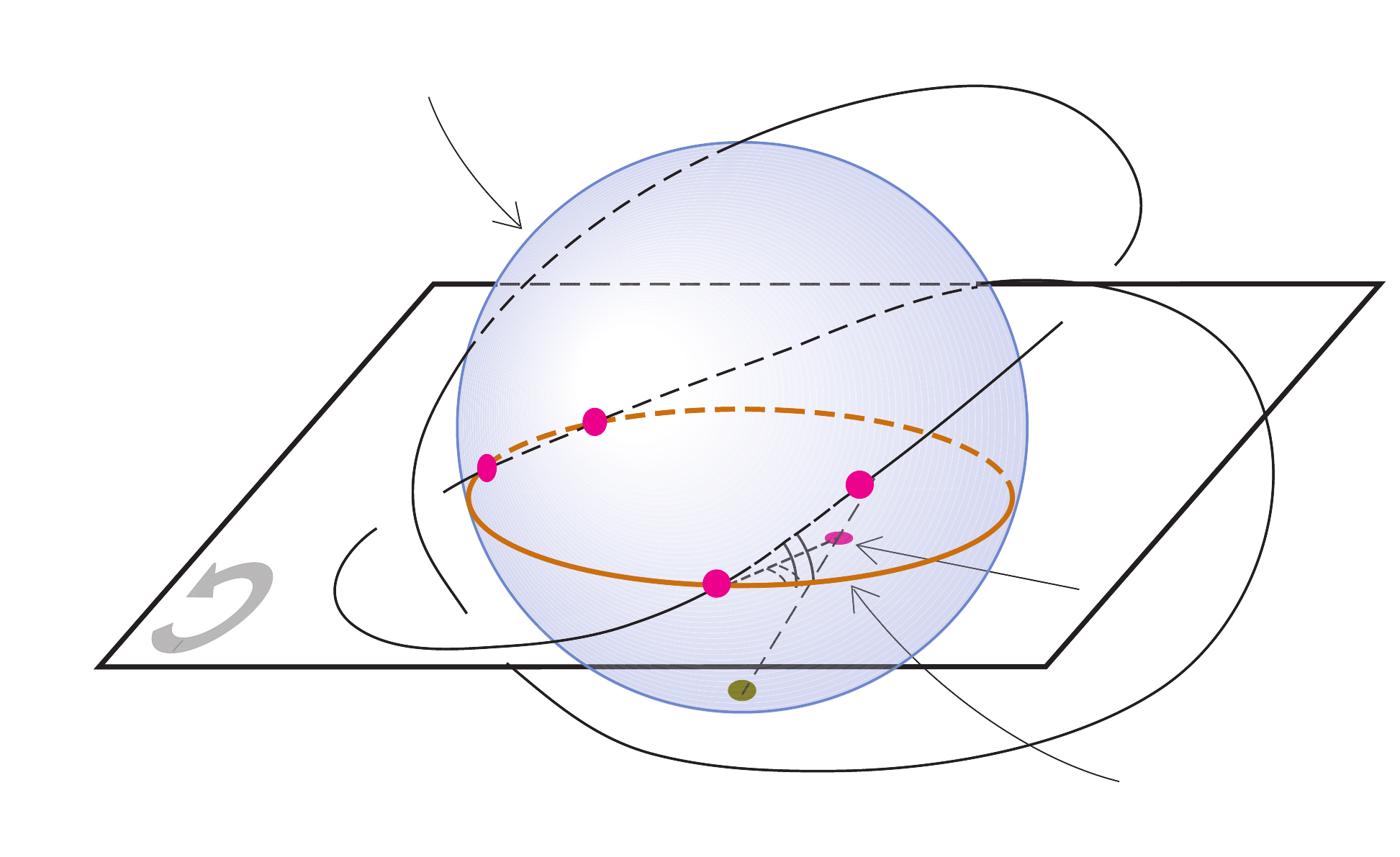}}
\caption{A stereographic projection}
\label{infcr_knot11c}
\end{figure}

Let $\varSigma=\varSigma_K(x,y)$ denote the $2$--sphere which is tangent to the knot $K$ at both $x$ and $y$. 
We call it a {\em bitangent sphere}. 
It can be considered the $2$--sphere $\varSigma(x, x+dx, y, y+dy)$ that
passes through four points $x, x+dx, y$, and $y+dy$ (\fullref{infcr_knot11c} left). 
It is generically determined uniquely unless these four points are cocircular, which is a condimension $2$ phenomenon. 
Identify $\Sigma_K(x,y)$ with the Riemann sphere
$\mathbb{C}\cup\{\infty\}$ through a stereographic projection $p$
(\fullref{infcr_knot11c} right). 
Then the four points $x, x+dx, y$, and $y+dy$ can be considered a quadruplet of complex numbers $\tilde x=p(x)$, $\tilde x+\widetilde{dx}=p (x+dx)$, $\tilde y=p(y)$, and $\tilde y+\widetilde{dy}=p(y+dy)$. 
Let $\Omega_K(x,y)$ be the cross ratio $(\tilde x+\widetilde{dx}, \tilde y ; \tilde x, \tilde y+\widetilde{dy})$:
$$
\Omega_K(x,y)
=\frac{(\tilde x+\widetilde{dx})-\tilde x}{(\tilde x+\widetilde{dx})
-(\tilde y+\widetilde{dy})}
\co \frac{\tilde y-\tilde x}{\tilde y-(\tilde y+\widetilde{dy})}
\sim\frac{\widetilde{dx}\widetilde{dy}}{(\tilde x-\tilde y)^2}. 
$$
To be precise, we need an orientation of $\varSigma$ to avoid the ambiguity of complex conjugacy of the infinitesimal cross ratio (see the Remark below). 

Then $\Omega_K(x,y)$ is independent of the choice of the stereographic projection $p$. 
Suppose we use another stereographic projection. 
Then we get another quadruplet of complex numbers, which can be obtained from the former by a linear fractional transformation. 
Since a linear fractional transformation does not change the cross ratio, we get the same value.

We call $\Omega_K(x,y)$ the {\em infinitesimal cross ratio} of the knot $K$. 
The real part of it has a pole of order $2$ at the diagonal $\Delta\subset K\times K$. 
It satisfies 
$$
(T\times T)^{\ast}\left(\Omega_{T(K)}(Tx,Ty)\right)=\Omega_K(x,y)
$$
for any M\"obius transformation $T$, where $T\times T$ is the diagonal action 
$$
T\times T\co K\times K\setminus\Delta\ni(x,y)\mapsto (Tx,Ty)\in T(K)\times T(K)\setminus\Delta.
$$ 
\begin{remark}\rm 
Let $\mathcal{S}$ be the set of quadruplets of ordered four points in
$S^3$ which are not cocircular.  We can define a continuous map from
$\mathcal{S}$ to the set of oriented spheres $\varLambda$.  (It is
given by a similar formula to \eqref{f_def_psi2} which shall be given
later.)  The composite with the cross ratio map gives a continuous map
from $\mathcal{S}$ to $\mathbb{C}\setminus\mathbb{R}$.  Since
$\mathcal{S}$ is connected its image is contained in one of the two
half-planes.  Our convention implies that the imaginary part of the
cross ratio of any ordered quadruplet of non-cocircular points in
$S^3$ is non-negative (the reader is referred to O'Hara \cite{OH5} for
details).
\end{remark}

\subsection{Conformal angles and cosine formula}
\begin{defn}\label{def_conf_angle} \rm (Doyle and Schramm)\qua
Let $C(x,x,y)$ be an oriented circle tangent to $K$ at $x$ which passes through $y$ whose orientation coincides with that of $K$ at $x$. 
Let $\theta$ be the angle from $C(x,x,y)$ to $C(y,y,x)$ at point $y$. 
We call it the {\em conformal angle} between $x$ and $y$ and denote it by $\theta=\theta_K(x,y)$.
\end{defn}

\begin{figure}[ht!]
\cl{\labellist\small
\pinlabel $K$ [b] at 262 323
\pinlabel $x$ [b] at 109 190
\pinlabel* $y$ [tl] at 168 114
\pinlabel $C(x,x,y)$ [t] <0pt, 3pt> at 73 46
\endlabellist
\includegraphics[height=4.4cm]{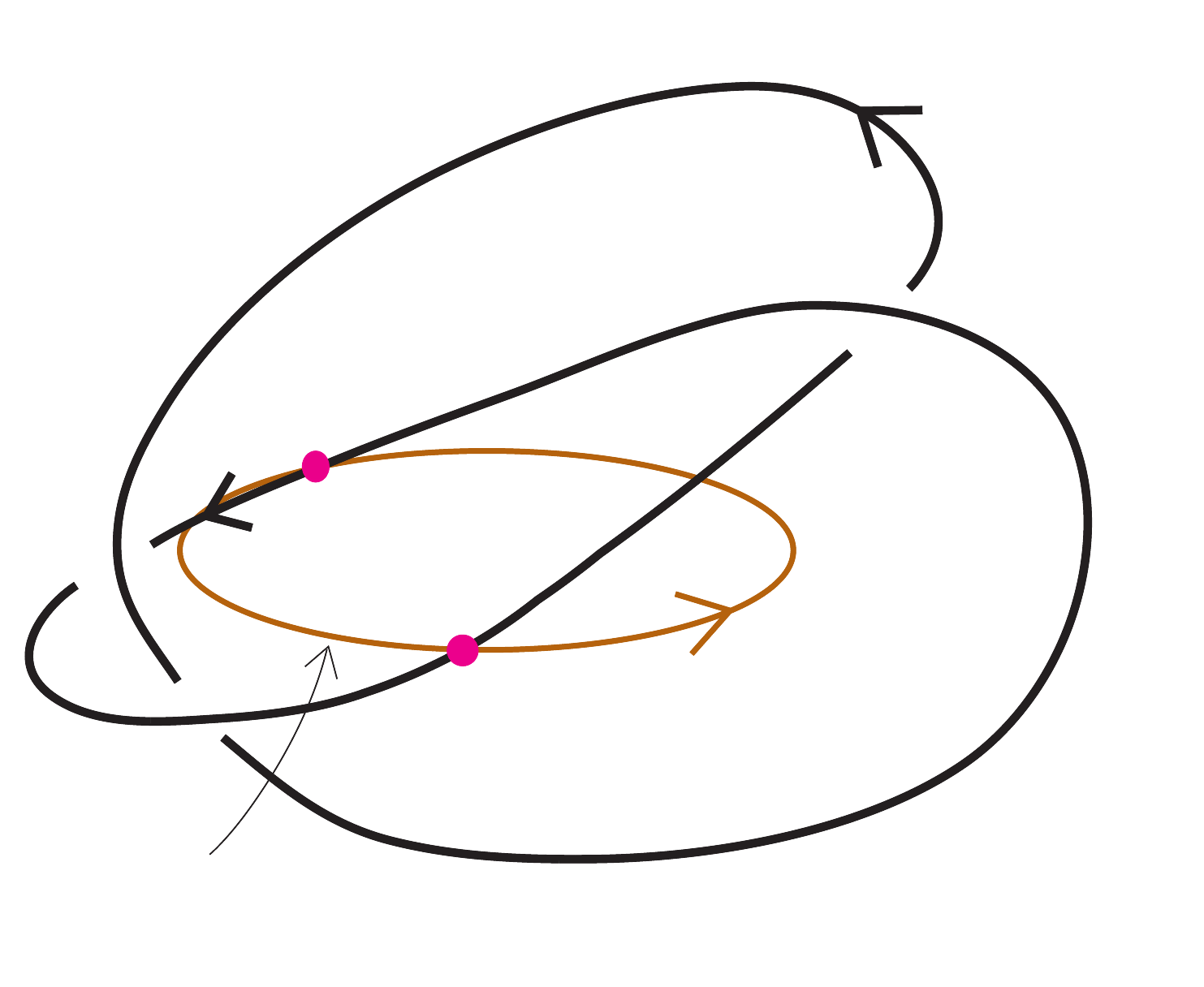}\qquad
\labellist\small
\pinlabel $K$ [b] at 262 323
\pinlabel $x$ [b] at 109 190
\pinlabel* $y$ [tl] at 168 114
\pinlabel $C(x,x,y)$ [t]  <0pt, 3pt> at 73 46
\pinlabel* $C(y,y,x)$ [tl]  <1.5pt, 0pt> at 293 186
\pinlabel* $\theta_K(x,y)$ [tl] at 257 113
\endlabellist
\includegraphics[height=4.4cm]{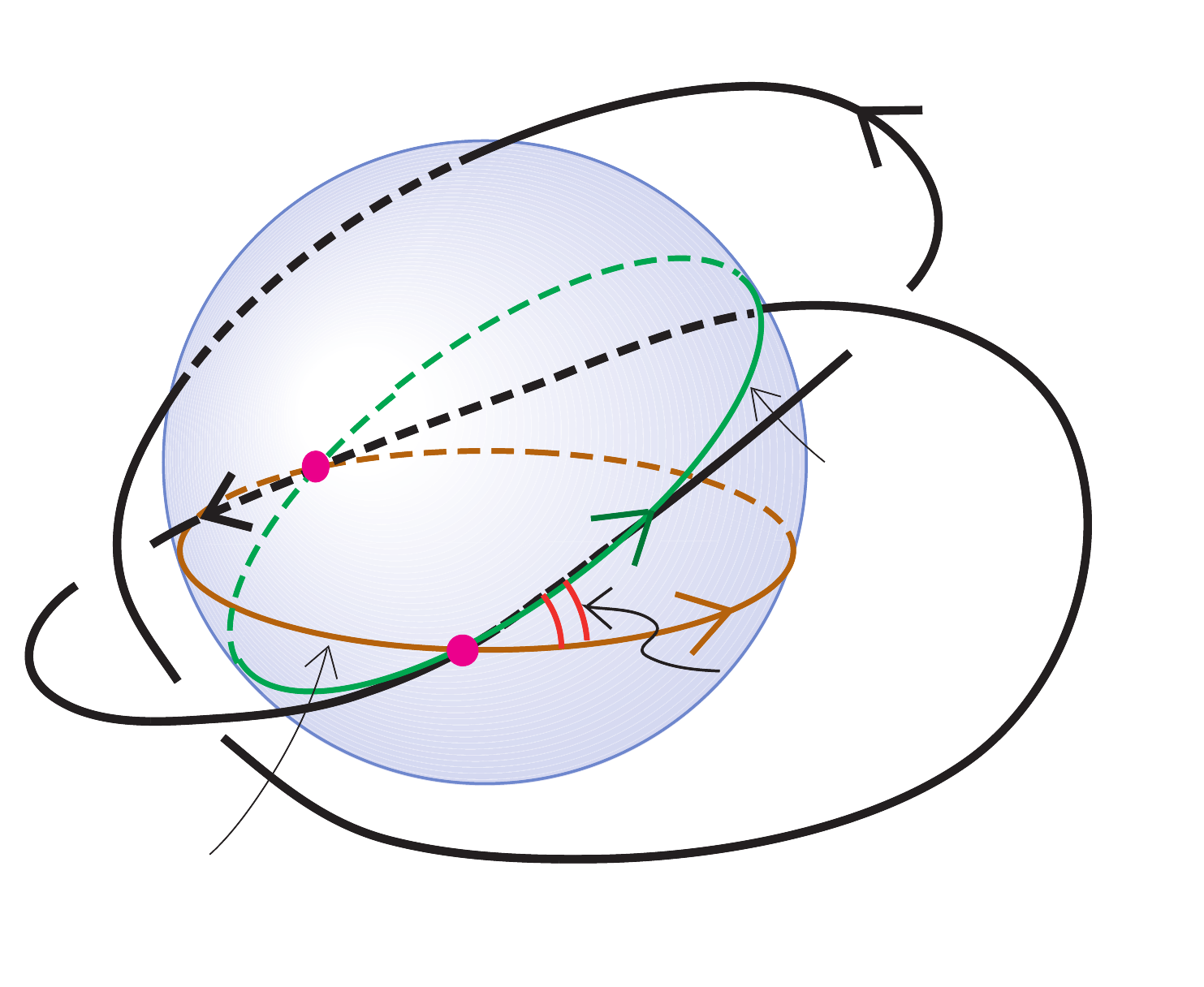}}
\caption{The conformal angle}
\label{conf_angl-c2-4}
\end{figure}

Generically $C(x,x,y)$ and $C(y,y,x)$ are different. 
The bitangnet sphere $\varSigma_K(x,y)$ is then the unique sphere that contains both $C(x,x,y)$ and $C(y,y,x)$. 
We assume that the sign of $\theta_K(x,y)$ is given with respect to the orientation of $\varSigma_K(x,y)$. 
Our convention of the orientation of bitangent spheres \cite{OH5} implies that the conformal angle always satisfies $0\le\theta_K(x,y)\le\pi$. 

\begin{proposition}\label{formula_inf_cr}
The absolute value of the infinitesimal cross ration $\Omega_K(x,y)$ is equal to $\displaystyle \frac{dxdy}{|x-y|^2}$ and the argument is equal to $\theta_K(x,y)$. 
Therefore we have 
$$\displaystyle \Omega_K(x,y)=e^{i\theta_K(x,y)}\frac{dxdy}{|x-y|^2}.$$ 
\end{proposition}

\begin{remark} 

(1)\qua The conformal angle is of the order of $|x-y|^2$ near the diagonal. 

(2)\qua The conformal angle behaves like an absolute value of a smooth function. 
Therefore, the imaginary part of the infinitesimal cross ratio may have singularity at $\{(x,y)\in K\times K\setminus\Delta\,|\,\theta_K(x,y)=0\}$. 
\end{remark}

Doyle and Schramm gave a {\em cosine formula} of $E^{(2)}_{\circ}(K)$
(Auckly and Sadun \cite{Au-Sa}, Kusner and Sullivan \cite{Ku-Su}): 
$$
E^{(2)}_{\circ}(K)=\iint_{\!\!K\times K\setminus\Delta}\frac{(1-\cos\theta_K(x,y))}{|x-y|^2}\,dxdy.
$$
This is another proof of the conformal invariance of $E_{\circ}^{(2)}$. 

\fullref{formula_inf_cr} and the cosine formula imply: 
\begin{proposition}{\rm (Langevin and O'Hara \cite{LO1})} The energy $E_{\circ}^{(2)}$ can be expressed in terms of the infinitesimal cross ratio $\Omega_{K}$ as 
$$\displaystyle 
E^{(2)}_{\circ}(K)=\iint_{K\times K\setminus\Delta} \left(|\Omega_{K}|-\R \Omega_{K}\right). 
$$
\end{proposition}

\subsection{$\R \Omega_K$ and the canonical symplectic form of $T^{\ast}S^3$}
We give the first interpretation of the real part of the infinitesimal cross ratio as the pull-back of the canonical symplectic form of the cotangent bundle $T^{\ast}S^3$. 

\begin{defn} \rm 
Let $T^{\ast}M$ be a cotangent bundle of an $m$--dimensional manifold $M$. 
Let $(q_1, \cdots, q_m, p_1, \cdots, p_m)$ be local coordinates of $T^{\ast}M$, where $(q_1, \cdots, q_m)$ are local coordinates of $M$ and $(p_1, \cdots, p_m)$ are local coordinates of fibers associated with the basis $\{dq_1, \cdots, dq_m\}$. 
The {\em canonical symplectic form $\omegasub{M}$ of the cotangent bundle} $T^{\ast}M$ is a globally defined non-vanishing $2$--form which can locally be expressed by 
$$\omegasub{M}=\sum dq_i\w dp_i.$$
\end{defn}

It is an exact form. 
In fact, there is a 1-form $\theta$ of $T^{\ast}M$ which can locally be expressed by $\theta=\sum p_idq_i$ that satisfies $\omegasub{M}=-d\theta$. 
This $\theta$ is called the tautological form of $T^{\ast}M$. 
It can be defined globally as follows. 
Let $T(T^{\ast}M)$ be a tangent bundle of $T^{\ast}M$. 
Then $\theta$ is given by 
$$(\theta(x,v))(w)=v(d\pi(w))\in\mathbb{R}, \hskip 0.5cm (x,v)\in T_x^{\ast}M,
 \hskip 0.3cm w\in T_{(x,v)}T^{\ast}M,$$
where $d\pi\co T_{(x,v)}T^{\ast}M\to T_xM$ is induced by the projection $\pi\co T^{\ast}M\to M$. 

The space $S^n\times S^n\setminus\Delta$ can be identified with the total space of the cotangent bundle $T^{\ast}S^n$ as follows. 
Assume $S^n\subset\mathbb{R}^{n+1}$. 
Let $\vect x\in S^n$. 
Let $\varPi_{\mbox{\scriptsize \boldmath$x$\normalsize}}=\left(\textrm{Span}\langle\vect x\rangle\right)^{\perp}$ be the $n$--plane in $\mathbb{R}^{n+1}$ through the origin which is orthogonal to $\vect x$, and $p_{\mbox{\scriptsize \boldmath$x$\normalsize}}\co S^n\setminus\{\vect x\}\to\varPi _{\mbox{\scriptsize \boldmath$x$\normalsize}}$ be a stereographic projection.
We identify $T_{\mbox{\scriptsize \boldmath$x$\normalsize}}S^n$ with $\varPi_{\mbox{\scriptsize \boldmath$x$\normalsize}}\cong\mathbb{R}^n$, and $T_{\mbox{\scriptsize \boldmath$x$\normalsize}}S^n$ with $T_{\mbox{\scriptsize \boldmath$x$\normalsize}}^{\ast}S^n$ by 
$$
T_{\mbox{\scriptsize \boldmath$x$\normalsize}}S^n\ni\vect u\mapsto \left(T_{\mbox{\scriptsize \boldmath$x$\normalsize}}S^n\ni\vect v\mapsto (\vect u,\vect v)\in\mathbb{R}\right)\in T_{\mbox{\scriptsize \boldmath$x$\normalsize}}^{\ast}S^n.
$$ 
Then the composition of identifications 
$$
\varphi_{\mbox{\scriptsize \boldmath$x$\normalsize}}\co S^n\setminus\{\vect x\}\spbmapright{\cong}{p_{\mbox{\tiny \boldmath$x$\normalsize}}}\varPi_{\mbox{\scriptsize \boldmath$x$\normalsize}}\spbmapright{\cong}{} T_{\mbox{\scriptsize \boldmath$x$\normalsize}}S^n \spbmapright{\cong}{} T_{\mbox{\scriptsize \boldmath$x$\normalsize}}^{\ast}S^n
$$
induces a canonical bijection $\varphi$: 
$$
S^n\times S^n\setminus\Delta
=\!\bigcup_{\mbox{\scriptsize \boldmath$x$\normalsize}\in S^n}\!\{\vect x\}\times (S^n\setminus\{\vect x\})
\ni (\vect x,\vect y)\stackrel{\varphi}{\mapsto} (\vect x,\varphi_{\mbox{\scriptsize \boldmath$x$\normalsize}}(\vect y))
\in \!\bigcup_{\mbox{\scriptsize \boldmath$x$\normalsize}\in S^n}\!T_{\mbox{\scriptsize \boldmath$x$\normalsize}}^{\ast}S^n
=T^{\ast}S^n.
$$ 
Let us write the pull-back $\varphi^{\ast}\omegasub{S^n}$ of the canonical symplectic form $\omegasub{S^n}$ of $T^{\ast}S^n$ by the same letter $\omegasub{S^n}$. 
\begin{thm}\label{thm_LO1}{\rm (Langevin and O'Hara \cite{LO1})} 

\begin{enumerate}
\item The $2$--form $\omegasub{S^n}$ on $S^n\times S^n\setminus\Delta$ is invariant under the diagonal action of a M\"obius transformation: $(T\times T)^{\ast}\omegasub{S^n}=\omegasub{S^n}$. 
\item {\rm (Folklore)}\qua Let $\displaystyle \omega_{\rm cr}\!=\!\frac{dw\!\w\! dz}{(w\!-\!z)^2}$ be a complex $2$--form on $\mathbb{C}\!\times\!\mathbb{C}\!\setminus\Delta$. 
It can be considered the cross ratio of $w,w+d w, z, z+d z$: 
$$
\displaystyle \frac{( w+{d w})- w}{( w+{d w})-( z+{d z})}
\co \frac{ z- w}{ z-( z+{d z})}
=\frac{{d w}{d z}}{( w- z)^2}\,. 
$$
Then $\R \omega_{\rm cr}=-\frac12\,\omegasub{S^2}$ through the identification of $S^2$ and $\mathbb{C}\cup\{\infty\}$ by a stereographic projection. 
\item 
The real part of the infinitesimal cross ratio can be expressed as the pull-back of the canonical symplectic form of the cotangent bundle $T^{\ast}S^3$ by the inclusion $\iota\co K\times K\setminus\Delta\hookrightarrow S^3\times S^3\setminus\Delta$: 
$$\displaystyle \R \Omega(x,y)=-\frac12\iota^{\ast}\omegasub{S^3}.$$
As a corollary, $\displaystyle \R \Omega(x,y)$ is an exact form. 
\end{enumerate}
\end{thm}
\subsection{Pseudo-Riemannian structure of the set of spheres}
We introduce some of the results from Langevin and O'Hara \cite{LO2}
in what follows.

We give the second interpretation of the real part of the infinitesimal cross ratio using the pseudo-Riemannian structure of the set of oriented $0$--spheres in $S^3$. 

Let $\mathcal{S}(q,n)$ denote the set of oriented $q$--spheres in $S^n$. 
As we saw in \fullref{sec_de_Sitter_sp_as_S}, when $n=3$ and $q=2=3-1$, $\mathcal{S}(2,3)$ can be identified with the de Sitter space $\varLambda$ in $\mathbb{R}^{5}_1$. 
The restriction of the indefinite metric of $\mathbb{R}^{5}_1$ to each tangent space of $\varLambda$ induces an indefinite non-degenerate quadratic form of index $1$. 
Let us consider the generalization to the cases with bigger codimensions. 
We assume $n-q\ge 2$ in this subsection. 

\begin{thm}\label{main_thm_LO2}{\rm \cite{LO2}}\qua
The dimension of $\mathcal{S}(q,n)$ is given by $(q+2)(n-q)$. 
There is a natural pseudo-Riemannian structure on $\mathcal{S}(q,n)$ of index $n-q$. 
Namely, each tangent space $T_p\mathcal{S}(q,n)$ admits an indefinite non-degenerate quadratic form $g$ such that $T_p\mathcal{S}(q,n)$ can be decomposed as the direct sum $T_p\mathcal{S}(q,n)\cong V_+\oplus V_-$ such that $\dim V_+=(q+1)(n-q)$, $\dim V_-=n-q$, and that the restriction of $g$ to $V_+$ {\rm (}or, to $V_-${\rm )} is positive definite {\rm (}or respectively, negative definite{\rm )}. 
\end{thm}
This indefinite non-degenerate quadratic form $g$ induces an indefinte pseudo-inner product. 

Just like in the case of $n=3$, $S^n$ can be realized in the Minkowski space $\mathbb{R}^{n+2}_1$ with the metric 
$$
\langle\vect x, \vect x\rangle=-x_0{}^2+x_1{}^2+\cdots+x_{n+1}{}^2
$$ 
as the set of lines through the origin in the light cone $V=\left\{\langle\vect v,\vect v\rangle=0\right\}$, and an oriented $q$--sphere $\varSigma$ in $S^n$ can be considered the intersection of $S^n$ and an oriented $(q+2)$--plane $\varPi_{\varSigma}$ through the origin. 
Therefore, $\mathcal{S}(q,n)$ can be identified with the Grassmann manifold 
$$
\widetilde{\textsl{Gr}}_-(q+2;\mathbb{R}^{n+2}_1)=\left\{\varPi\subset\mathbb{R}^{n+2}_1\left|\!\begin{array}{l}
\mbox{oriented $(q+2)$--plane through $\vect 0$}\\
\mbox{$\varPi$ intersects the light cone transversally$^{\ast}$}
\end{array}\!\!\right\}\right..
$$
($\ast$ The second condition above is equivalent to say that $\varPi$ is again a Minkowski space, ie, the restriction of $\langle \,,\,\rangle$ to $\varPi$ is a non-degenerate indefinite quadratic form of index $1$.) 
It is a homogeneous space 
$$
\widetilde{\textsl{Gr}}_-(q+2;\mathbb{R}^{n+2}_1)\cong SO(n+1,1)/SO(n-q)\times SO(q+1,1), 
$$ and Theorem follows from Kobayashi and Yoshino \cite[Propopsition
3.2.6]{KY}.

Let us introduce more constructive explanation which is useful in the study of conformal geometry. 

Let $\varPi$ be an oriented $(q+2)$--dimensional vector subspace in
$\mathbb{R}^{n+2}_1$, and let $\{\vect x_1, \cdots,$\break $\vect x_{q+2}\}$
be an ordered basis of $\varPi$ which gives the orientation of
$\varPi$.  Let $M$ be a $(q+2)\times (n+2)$--matrix given by
$$
M=\left(\begin{array}{c}
\vect x_1\\
\vdots\\
\vect x_{q+2}
\end{array}\right)
=\left(\begin{array}{cccc}
x_{1\,0} & x_{1\,1} & \cdots & x_{1\,n+1}\\
\vdots & \vdots & \ddots & \vdots\\
x_{q+2\,0} & x_{q+2\,1} & \cdots & x_{q+2\,n+1}
\end{array}\right).
$$
Let $I=(i_1, \cdots, i_{q+2})$ be a multi-index $(0\le i_k\le n+1)$. 
Define $p_I=p_{i_1\cdots i_{q+2}}$ by 
\begin{eqnarray}\label{f_def_plucker_coordinates}
p_{i_1\cdots i_{q+2}}=
\left|\begin{array}{cccc}
\,x_{1\,i_1}\,&\,x_{1\,i_2}\,&\cdots&\,x_{1\,i_{q+2}}\,\\
x_{2\,i_1}&x_{2\,i_2}&\cdots&x_{2\,i_{q+2}}\\
\vdots&\vdots&\ddots&\vdots \\
x_{(q+2)\,i_1}&x_{(q+2)\,i_2}&\cdots&x_{(q+2)\,i_{q+2}}\\
\end{array}\right|.
\end{eqnarray}
Then $p_{i_1\cdots i_{q+2}}$ is alternating in the suffixes $i_k$. 
The exterior product of $\vect x_1, \cdots, \vect x_{q+2}$ is given by
$$\vect x_1\w\cdots\w\vect x_{q+2}
=\sum_{0\le i_1<\cdots <i_{q+2}\le n+1}p_{i_1\cdots i_{q+2}}\,\vect e_{i_1}\w\cdots\w\vect e_{i_{q+2}}\in\stackrel{q+2}{\bigwedge}\mathbb{R}^{n+2}_1.$$
Let $\displaystyle N={{n+2}\choose{q+2}}$. 
We identify $\displaystyle \stackrel{q+2}{\mbox{\Large$\wedge$}}\mathbb{R}^{n+2}_1$ with $\mathbb{R}^N$ by expressing $\vect x_1\w\cdots\w\vect x_{q+2}\in \displaystyle \stackrel{q+2}{\mbox{\Large$\wedge$}}\mathbb{R}^{n+2}_1$ by $(\cdots,p_{i_1\cdots i_{q+2}},\cdots)\in\mathbb{R}^N$. 

Let $[\varPi]$ denote an unoriented $(q+2)$--space which is obtained from $\varPi$ by forgetting its orientation. 
Then it can be identified by the homogeneous coordinates $[\cdots,p_{i_1\cdots i_{q+2}},\cdots]\in\mathbb{R}P^{N-1}$. 
They are called the {\em Pl\"ucker coordinates} or {\em Grassmann coordinates}. 
They do not depend on the choice of $(q+2)$ linearly independent vectors which span $[\varPi]$. 
Let ${\textsl{Gr}}\,(q+2,n+2)$ be the Grassmann manifold of the set of all $(q+2)$--dimensional vector subspaces in $\mathbb{R}^{n+2}$. 
The mapping 
$$
{\textsl{Gr}}\,(q+2,n+2)\ni[\varPi]\mapsto [\cdots,p_{i_1\cdots i_{q+2}},\cdots]\in\mathbb{R}P^{N-1}
$$
is called the {\em Grassmann mapping}. 

The Pl\"ucker coordinates $p_{i_1\cdots i_{q+2}}$ are not independent. 
They satisfy the {\em Pl\"ucker relations}: 
\begin{equation}\label{plucker_relations}
\sum_{k=1}^{q+3}(-1)^{k}p_{i_1\cdots i_{q+1}j_k}p_{j_1\cdots \widehat{j_k}\cdots j_{q+3}}=0,
\end{equation}
where $\widehat{j_k}$ indicates that the index $j_k$ is being removed. 
(We remark that the digits $i_1, \cdots , i_{q+1}, j_k$ in the multi-index above are not necessarily ordered according to their sizes. ) 
All the Pl\"ucker relations are not necessarily independent. 

The pseudo-Riemannian structure of $\displaystyle \stackrel{q+2}{\mbox{\Large$\wedge$}}\mathbb{R}^{n+2}_1$ is given by 
$$
\langle \vect e_{i_1}\w\cdots\w\vect e_{i_{q+2}}, \,\vect e_{j_1}\w\cdots\w\vect e_{j_{q+2}}\rangle
=-\left|\!
\begin{array}{ccc}
\langle \vect e_{i_1}, \vect e_{j_1}\rangle & \cdots & \langle \vect e_{i_1}, \vect e_{j_{q+2}}\rangle\\
\vdots & \ddots & \vdots \\
\langle \vect e_{i_{q+2}}, \vect e_{j_1}\rangle & \cdots & \langle \vect e_{i_{q+2}}, \vect e_{j_{q+2}}\rangle\\
\end{array}\!
\right|,
$$
which can be obtained by generalizing a formula in the case of codimension $1$,\break $\displaystyle \stackrel{4}{\mbox{\Large$\wedge$}}\mathbb{R}^{5}_1\cong\mathbb{R}^{5}_1$. 
Therefore, 
$$\{\vect e_{i_1}\w\cdots\w\vect e_{i_{q+2}}\}_{0\le i_1<\cdots <i_{q+2}\le n+1}$$
can serve as a pseudoorthonormal basis of $\displaystyle \stackrel{q+2}{\mbox{\Large$\wedge$}}\mathbb{R}^{n+2}_1$ which satisfies  
$$
\langle \vect e_{i_1}\w\cdots\w\vect e_{i_{q+2}}, \,\vect e_{i_1}\w\cdots\w\vect e_{i_{q+2}}\rangle=\left\{
\begin{array}{lcl}
-1 & \textrm{ if } & i_1\ge 1,\\[1mm]
+1 & \textrm{ if } & i_1=0.
\end{array}
\right.
$$
It follows that if $\vect v=(\cdots,p_{i_1\cdots i_{q+2}},\cdots)\in\displaystyle \stackrel{q+2}{\mbox{\Large$\wedge$}}\mathbb{R}^{n+2}_1\cong\mathbb{R}^N$ then 
\begin{equation}\label{f_pR_str}
\langle \vect v, \vect v\rangle=-\sum_{1\le i_1<\cdots<i_{q+2}}p_{i_1\cdots i_{q+2}}{}^2+\sum_{i_1=0<i_2<\cdots<i_{q+2}}p_{0i_2\cdots i_{q+2}}{}^2. 
\end{equation}
Put $N_1={{n+1}\choose{q+2}}$ and $N_2={{n+1}\choose{q+1}}$. 
Then $\displaystyle \stackrel{q+2}{\mbox{\Large$\wedge$}}\mathbb{R}^{n+2}_1\cong\mathbb{R}^N$ can be decomposed to a direct sum $\mathbb{R}^{N_1}_-\oplus\mathbb{R}^{N_2}_+$, where the restriction of $\langle\,,\,\rangle$ to $\mathbb{R}^{N_1}_-$ (or $\mathbb{R}^{N_2}_+$) is negative (or respectively, positive) definite. 
We denote $\mathbb{R}^N$ with the metric $\langle\,,\,\rangle$ given by \eqref{f_pR_str} by $\mathbb{R}^N_{N_1}$. 

\begin{thm}{\rm \cite{LO2}}\qua 
Let $N={{n+2}\choose{q+2}}$ and $N_1={{n+1}\choose{q+2}}$ as before. 

\medskip
{\rm (1)}\qua Let $\varPi=\textsl{Span}\langle\vect x_1, \cdots, \vect x_{q+2}\rangle$ be an oriented $(q+2)$--dimensional vector subspace in $\mathbb{R}^{n+2}_1$ spanned by $\vect x_1, \cdots, \vect x_{q+2}$. 
Put $\vect p=\vect x_1\w\cdots\w\vect x_{q+2}\in\mathbb{R}^N_{N_1}$. 
Then $\varPi$ intersects the light cone $V$ transversally if and only if $\langle \vect{p}, \vect{p}\rangle>0$. 

\medskip{\rm (2)}\qua Let $S^{N-1}_{N_1}$ be the unit pseudosphere: 
$$
S^{N-1}_{N_1}=\left.\left\{\vect v=(\cdots,p_{i_1\cdots i_{q+2}},\cdots)\in\mathbb{R}^N_{N_1}\right|\langle \vect v, \vect v\rangle=1\right\} 
$$
and $\widetilde Q_P(q+2;\mathbb{R}^{n+2}_1)$ be the quadric satisfying the Pl\"ucker relations: 
$$
\widetilde Q_P(q+2;\mathbb{R}^{n+2}_1)=\left\{(\cdots,p_{i_1\cdots i_{q+2}},\cdots)\,\left|\,\sum_{k=1}^{q+3}(-1)^{k}p_{i_1\cdots i_{q+1}j_k}p_{j_1\cdots \widehat{j_k}\cdots j_{q+3}}=0\right\}\right..
$$
Then the set $\mathcal{S}(q,n)$ of oriented $q$--dimensional spheres in $S^n$ can be identified with the intersection of $S^{N-1}_{N_1}$ and $\widetilde Q_P(q+2;\mathbb{R}^{n+2}_1)$: 
$$
\mathcal{S}(q,n)\cong S^{N-1}_{N_1}\cap\widetilde Q_P(q+2;\mathbb{R}^{n+2}_1)\subset \mathbb{R}^N_{N_1}.
$$
Let us denote the right hand side by $\Theta(q,n)$. 

\medskip{\rm (3)}\qua Let $\varSigma(\vect x_1, \cdots, \vect x_{q+2})$ denote an oriented $q$--sphere $\varSigma$ which is given as the intersection of $S^n$ and an oriented vector subspace $\textsl{Span}\langle\vect x_1, \cdots, \vect x_{q+2}\rangle$. 
Then the bijection $\psi_G\co \mathcal{S}(q,n)\to\Theta(q,n)$ is given by
\begin{equation}\label{f_def_psi2}
\psi_G(\varSigma(\vect x_1, \cdots, \vect x_{q+2}))=\frac{\vect x_1\w\cdots\w\vect x_{q+2}}{\,\sqrt{\langle\vect x_1\w\cdots\w\vect x_{q+2},\vect x_1\w\cdots\w\vect x_{q+2}\rangle\,}\,}\,. 
\end{equation}
\end{thm}
We show that a M\"obius transformation of $S^n$ induces a pseudoorthogonal transformation of $\Theta(q,n)\subset\displaystyle \stackrel{q+2}{\mbox{\Large$\wedge$}}\mathbb{R}^{n+2}_1$. 
Let $O(N_2,N_1)$ denote the pseudoorthogonal group. 
\begin{defn} \rm Let $N={{n+2}\choose{q+2}}$ as before. 
Define 
$$
\varPsi_{q,n}\co M_{n+2}(\mathbb{R})\ni A=(a_{ij})\mapsto \widetilde{A}=(\tilde a_{IJ})\in M_N(\mathbb{R}),
$$
where $I=(i_1\cdots i_{q+2})$ and $J=(j_1\cdots j_{q+2})$ are multi-indices, 
and $\tilde a_{IJ}$ is given by 
$$
\tilde a_{IJ}=
\left|\!
\begin{array}{ccc}
\,a_{i_{1}j_{1}}\,& \cdots & \,a_{i_{1}j_{q+2}}\,\\
\vdots & \ddots & \vdots \\
a_{i_{q+2}j_{1}}& \cdots & a_{i_{q+2}j_{q+2}}
\end{array}\!
\right|. 
$$
\end{defn}
\begin{proposition}{\rm \cite{LO2}}\label{prop_tilde_A}\qua
Let $N={{n+2}\choose{q+2}}$, $N_1={{n+1}\choose{q+2}}$ and $N_2={{n+1}\choose{q+1}}$ as before. 
Let $A=(a_{ij})\in M_{n+2}(\mathbb{R})$. 
\begin{enumerate}
\item A matrix $\widetilde{A}\in M_N(\mathbb{R})$ satisfies 
$$
(A\vect x_1)\w \cdots \w(A\vect x_{q+2})=\widetilde{A}\,(\vect x_1\w \cdots\w \vect x_{q+2})\, (\forall \vect x_1, \cdots, \vect x_{q+2}\in\mathbb{R}^{n+2}_1) 
$$
if and only if $\widetilde{A}=\varPsi(A)$. 
\item If $A\in O(n+1,1)$ then $\varPsi(A)\in O(N_2,N_1)$. 
\item The restriction of $\varPsi_{q,n}$ to $O(n+1,1)$, which shall also be denoted by $\varPsi_{q,n}$, 
$$\varPsi_{q,n}\co O(n+1,1)\ni A\mapsto \widetilde{A}\in O(N_2,N_1)$$
is a homomorphism. 
\item Let $\psi_G\co \mathcal{S}(q,n)\to\Theta(q,n)$ be the bijection given in the above {\rm Theorem}. 
Then, 
$$\psi_G(A\cdot\varSigma)=\varPsi_{q,n}(A)\psi_G(\varSigma)$$
for $\varSigma\in\mathcal{S}(q,n)$ and  $A\in O(n+1,1)$. 
\end{enumerate}
\end{proposition}

\subsection{$\R \Omega$ as a signed area form}
We give the second interpretation of the real part of the infinitesimal cross ratio. 

\begin{figure}[ht!]
\begin{center}
\includegraphics[width=.66\linewidth]{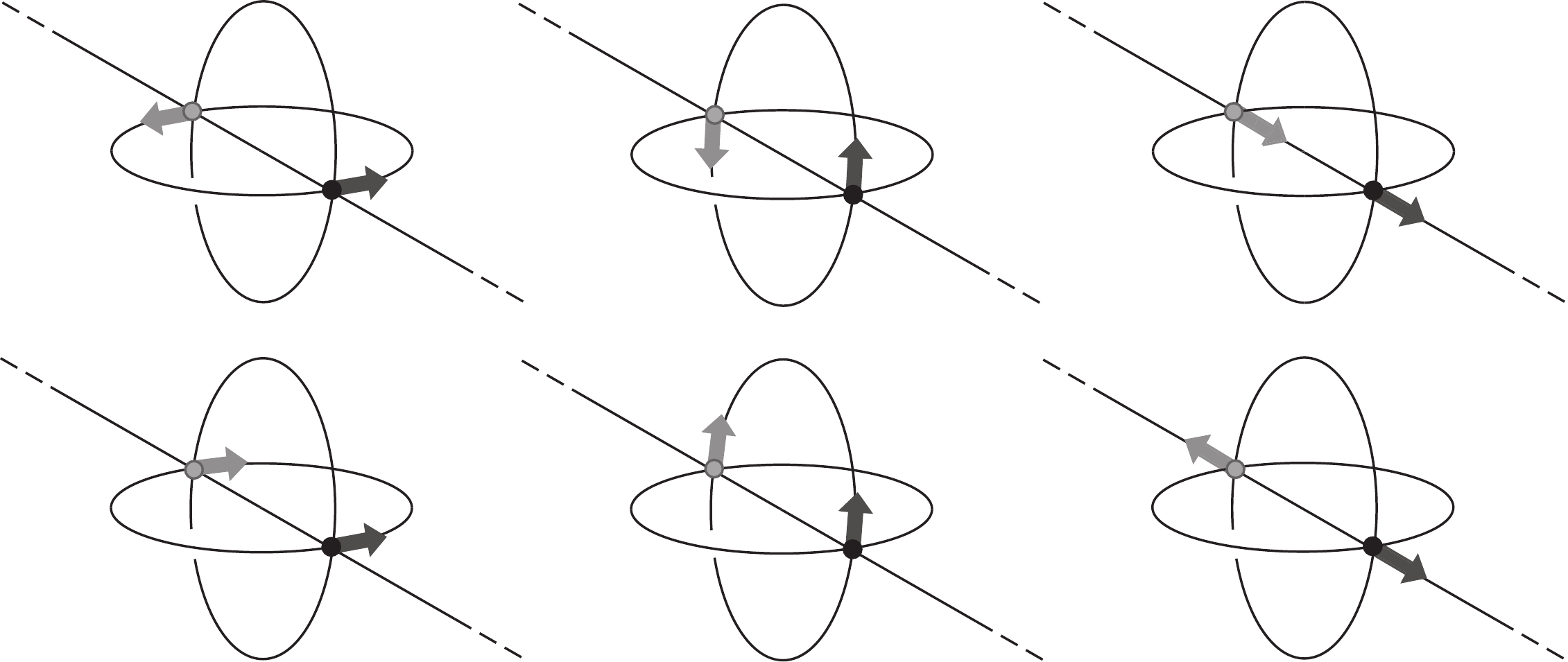}
\caption{$3$ spacelike pencils (above) and $3$ {timelike} pencils (below)}
\label{orthonormal_basis_S^3-2}
\end{center}
\end{figure}

When $q=0$ and $n=3$ the set $\mathcal{S}(0,3)$ of oriented $0$--spheres in $S^3$ is a subspace of $\mathbb{R}^{10}_6$ since $N={{5}\choose{2}}=10$ and $N_1={{4}\choose{2}}=6$. 
At the same time, it is identical with $S^3\times S^3\setminus\Delta$. 
It admits the pseudo-Riemannian structure of index $3$ (\fullref{main_thm_LO2}). 
The pseudoorthonormal basis can be given by mutually pseudoorthogonal {\sl pencils}, as is illustrated in \fullref{orthonormal_basis_S^3-2} (which is a picture in $\mathbb{R}^3$ obtained through a stereographic projection). 

Let $(x,y)$ be a pair of distinct points of a knot $K$. 
Then it can be considered a point in $\mathcal{S}(0,3)\cong\Theta(0,3)$. 
Let it be denoted by $\vect s(x,y)$. 
Namely, $\vect s$ induces a map
$$\vect s\co K\times K\setminus\Delta\hookrightarrow S^3\times S^3\setminus\Delta\stackrel{\cong}{\longrightarrow}\mathcal{S}(0,3)\cong\Theta(0,3)\subset\mathbb{R}^{10}_6.$$ 
The image $\vect s(K\times K\setminus\Delta)$ is a surface in $\Theta(0,3)$. 
Its area element is given by 
$$
dv=\sqrt{\,\left|\!\begin{array}{cc}
\langle \vect s_x, \vect s_x \rangle  &  \langle \vect s_x, \vect s_y \rangle \\
\langle \vect s_y, \vect s_x \rangle & \langle \vect s_y, \vect s_y \rangle 
\end{array}\!\right|\,}\,dxdy\,,
$$
where $\vect s_x$ and $\vect s_y$ denote $\displaystyle \frac{\partial \vect s}{\partial x}(x,y)$ and $\displaystyle \frac{\partial \vect s}{\partial y}(x,y)$ in $T_{\vect s(x,y)}\Theta(0,3)$. 
It turns out that $\langle \vect s_x, \vect s_x \rangle=\langle \vect s_y, \vect s_y \rangle=0$. 
Therefore 
$$dv=\sqrt{-\langle \vect s_x, \vect s_y \rangle^2}\,dxdy.$$ 
\begin{defn} \rm 
Define a {\em signed area form} $\alpha$ of the surface $\vect s(K\times K\setminus\Delta)$ by 
$$
\alpha=\langle \vect s_x, \vect s_y \rangle \, dx\wedge dy. 
$$
\end{defn}
\begin{thm}{\rm \cite{LO2}}\qua
The real part of the infinitesimal cross ratio is equal to the half of the signed area form of the surface $\vect s(K\times K\setminus\triangle)$ with respect to the pseudo-Riemannian structure of $\mathcal{S}(0,3)$: 
$$
\R\Omega_K(x,y)=\frac{1}{\,2\,}\,\langle \vect s_x, \vect s_y \rangle \, dx\wedge dy.
$$
\end{thm}

Let $\gamma_1\cup\gamma_2$ be a $2$--component link. 
The infinitesimal cross ratio $\Omega(x,y)$ $(x\in\gamma_1, y\in\gamma_2)$ can be defined in the same way. 
The above Theorem and \fullref{thm_LO1} imply that the signed area form $\alpha=2\,\R\Omega$ of the surface $\vect s(\gamma_1\times\gamma_2)\subset\Theta(0,3)$ is an exact form. 
Therefore, Stokes' theorem implies that the signed area of $\vect s(\gamma_1\times\gamma_2)$ vanishes: 
$$
\int_{\gamma_1\times\gamma_2}\alpha=\int_{x\in\gamma_1}\int_{y\in\gamma_2}\langle \vect s_x, \vect s_y \rangle \, dx\wedge dy=0. 
$$

\subsection{The imaginary part of the infinitesimal cross ratio}
Unlike the real part, the imaginary part $\I \Omega$ of the infinitesimal cross ratio does not have a nice global interpretation. 
It cannot be expressed as a pull-back of a globally defined $2$--form. 
(We cannot generalize the imaginary part of $\displaystyle \omega_{\rm cr}\!=\!\frac{dw\!\w\! dz}{(w\!-\!z)^2}$ to $S^n\times S^n\setminus\Delta$ for $n\ge 3$.) 
It might be singular at $(x,y)\in K\times K\setminus\triangle$ where 
the conformal angle $\theta_K(x,y)$ vanishes. 

The imaginary part $\I \Omega$ of the infinitesimal cross ratio can be considered a local transversal area element of geodesics in $\mathbb{H}^4$ joining pairs of points on the knot $K$. 
To be precise, let $S^3\cong\partial \mathbb{H}^4$, $l(x,y)$ be a geodesic in $\mathbb{H}^4$ joining a pair of points $x$ and $y$ on $K$, $\varPi_0$ be any totally geodesic $3$--space of $\mathbb{H}^4$ which is perpendicular to $l(x_0,y_0)$, and $S(x,y)=l(x,y)\cap\varPi_0$ be a surface in $\varPi_0$. 
Then $\I \Omega(x_0,y_0)$ is equal to the quater of the area element of $S(x,y)$ at $(x_0,y_0)$ \cite{LO2}.

\bibliographystyle{gtart}
\bibliography{link}

\end{document}